\documentclass[11pt]{elsarticle}
\usepackage{epsfig}
\usepackage{color}
\usepackage{amsmath}
\usepackage{amssymb}  %
\usepackage{comment}
\usepackage{pgfplots}
\usepackage{mathdots}
\usepackage{tikz}
\usepackage{url}
\usepackage{multirow}
\usepackage{appendix}

\usepackage{footnote}
\makesavenoteenv{tabular}
\makesavenoteenv{table}

\newfont{\Bb}{msbm10 scaled\magstep0}

\newtheorem{remark}{Remark}[section]

\usepackage{algorithm}
\usepackage{algorithmic}

\counterwithin{figure}{section}
\counterwithin{table}{section}

\begin{document}

\begin{frontmatter}

\title{Surrogate modeling for convection-dominated parametric problems based on error learning}

\author{Lihong Feng\corref{cor1}}
\ead{feng@mpi-magdeburg.mpg.de}
\author{Ali Hassan}
\author{Shuwen Sun}

\affiliation{organization={Max Planck Institute for Dynamics of Complex Technical Systems},
             city={Magdeburg},
             postcode={39016},
             country={Germany}}

\cortext[cor1]{Corresponding author}




\begin{abstract}
Convection-dominated problems are known for their slow Kolmogorov $n$-width decays and are challenging for both projection-based linear model order reduction (MOR) methods and data-driven MOR methods. In this work, we propose a hybrid surrogate modeling approach and a non-intrusive variant that overcome some drawbacks of linear MOR methods. The proposed hybrid surrogate model is a projection-based reduced-order model (ROM), corrected by the error learned from a deep neural network. With the aid of deep learning, the model component of the surrogate model can be kept in a small reduced dimension. The neural network component and the model component of the hybrid surrogate are sequentially but separately built during the offline stage. At the online stage, they are easily coupled to output the solution predictions at desired parameter values and time instances. Due to the intrusive nature of the hybrid-ROM, the numerically discretized operators of the original model must be available. For problems solved using black-box solvers, where the details of the numerical discretization are not accessible, we further propose a non-intrusive variant of the hybrid surrogate. Compared to the existing MOR methods with nonlinear manifolds, which couple the neural network within the ROM, the proposed hybrid ROM is more easily built and is also easily assembled for online prediction. Thanks to the physics-based input, the neural network component of the hybrid ROM accurately learns the error within a few hundreds to thousands of epochs. In contrast to the surrogate modeling approaches purely based on deep-learning, the proposed non-intrusive variant has a lighter neural network structure with much fewer parameters to be learned. We test the proposed methods on two nonlinear convection parametric problems. The first is the 1D inviscid Burgers' equation with one parameter, and the second is the 2D inviscid Burgers' equation with two parameters. Since both methods are based on error correction, their online predictions exhibit higher accuracy yet with largely reduced prediction time, compared to state-of-the-art methods.
\end{abstract}

\end{frontmatter}

\section{Introduction}
In the history of model order reduction for parametric time-dependent problems, many methods have primarily focused on problems whose solution can be well approximated by a linear subspace~\cite{morBenGW15}. For convection-dominated problems, traditional projection based MOR methods using linear subspace~\cite{morBenGW15} or recent neural network based surrogate modeling methods~\cite{morHesU18} relying only on proper orthogonal decomposition (POD) loose accuracy, especially when the linear subspace dimensions are kept small. The solution manifolds of the convection-dominated problems are usually with slow Kolmogorov $n$-width decay, which are challenging for the aforementioned linear-subspace-based methods. In the past years, quite a few MOR methods are proposed to tackle the challenges of convection-dominated problems.  Ideas of space/snapshot transforms are proposed~\cite{morBarTT24, morTadZ21, morRimPM23, morGruS24, morCagMS19, morOhlR13, morNaiB19, morEhrLMetal20, morSalDM21, morKhaPetal26} to capture the mass transport properties of the solutions and derive ROMs with small sizes.  While these methods show efficiency for some problems, most of them are limited to special situations, e.g., steady-state problems~\cite{morNaiB19, morBarTT24}, 1D problems~\cite{morEhrLMetal20, morTadZ21, morRimPM23, morCagMS19}, or problems with a single parameter~\cite{morKhaPetal26}. Some methods~\cite{morGruS24, morSalDM21} present results for 2D problems with yet  insufficient accuracy, while others lack practical implementations~\cite{morOhlR13}. 

Local adaptive methods are also proposed to update basis online~\cite{morPeh20}, but with higher demands on online computations, compared to the traditional MOR methods. Instead of space/snapshot transform, nonlinear compression using autoencoders are proposed to get projection-based ROMs for convection-dominated problems~\cite{morLeeC20, morKimCWZ22}. Non-intrusive methods using nonlinear compression via convolutional autoencoders (CAEs) have also been proposed to deliver fast neural-network-based surrogate models to deliver fast neural-network-based surrogate models~\cite{morFreetal21, morFreM22}, demonstrating superior performance over traditional POD-based ROMs. In parallel, data-driven closure modeling approaches~\cite{morAhmPSetal21, morMouKSetal21, morManTAetal25,  morKocRRetal25} claim their advantages over traditional MOR methods using linear subspaces. However, the contributions from the closure terms are sometimes marginal. This behavior may stem from the ill-conditioned least-squares (LS) problems employed in computing the closure terms. Note that the ROMs with closures~\cite{morAhmPSetal21, morMouKSetal21, morManTAetal25,morKocRRetal25} share similar formulations with the hybrid models proposed in~\cite{morChiCChetal20}. 

Projection-based methods using nonlinear manifold augmentation construct ROMs using combination of linear subspaces and nonlinear manifolds~\cite{morBarFM2023, morParTHetal26, morBarC22,  morAghBMetal26}. These methods, especially the methods using manifolds with general nonlinear structures~\cite{morBarFM2023,morParTHetal26} show highly improved accuracy over traditional linear MOR methods based on POD for pure convection problems and turbulence problems. One observation is that these methods require involved offline and online implementations, especially when hyper-reduction is unavoidable to construct the final nonlinear-manifold-augmented ROMs and when implicit time-integration methods are needed for accurately solving the ROMs. Nonlinear manifold based operator inference methods~\cite{morGeeBWetal24} use nonlinear manifolds represented by quadratic/polynomial basis to augment the linear subspace representation of the solution. These methods can be seen as the non-intrusive counterparts of the projection-based methods with quadratic manifold augmentation~\cite{morBarC22, morAghBMetal26}. The ROM operators are learned by formulating a LS problem. Efficiently solving the LS problem becomes the key issue, which is nevertheless often ill-conditioned and requires careful regularization techniques. One limitation of the methods based on nonlinear-manifold augmentation is that the accuracy of the ROMs is limited by the dimension of the augmented linear subspace, i.e., the constructed ROMs are no more accurate than the ROMs constructed by the linear subspace with the augmented dimension.

We propose a hybrid MOR approach that involves only standard MOR methods, e.g, the reduced-basis method (RBM), and a neural network. The structure of the hybrid ROM is simple, $\tilde {\bf u}(t,\mu)=V {\bf z}(t,\mu)+\tilde {\bf e}({\bf z}(t, \mu))$, where ${\bf z}(t,\mu)$ is the solution to a RB-ROM, and $\tilde{\bf e}({\bf z}(t,\mu))$ is the error learned by a neural network. Matrix $V\in \mathbb R^{n\times r}, r\ll n$ is the projection matrix computed from the RBM, where $n$ is the dimension of the full-order model (FOM). The proposed method breaks the limitation of the nonlinear manifold augmented MOR methods by directly learning the solution error of the model component of the hybrid ROM. For problems where intrusive MOR cannot be applied, we further propose a non-intrusive variant that is purely data-driven. Instead of obtaining the latent state from the projection-based ROM, we compute the latent state from POD of the snapshots. A second neural network is established to map the parameter domain and the time domain to the latent space, which is trained by optimizing the error between the POD latent states and the NN-learned latent states.

The key advantages of the proposed surrogate modeling approaches over the existing MOR methods for convection-dominated problems are the simple model structure leading to much simpler hyper-reduction procedure, simple implementation procedure, simple coupling, physics-informed decoder, and applicability to parametric time-dependent nonlinear problems with multiple parameters. Numerical tests are done on the 1D inviscid Burgers' equation with parametrized initial condition and the 2D inviscid Burgers' equation with two parameters in the source term and on the boundary, respectively. The results show the effectiveness of the proposed method as compared to the state-of-the-art approaches. 

In the following sections, we first introduce the problem setting and the RB-ROM for parametric time-dependent problems in Section~\ref{sec:POD-greedy}. At the end of this section, we also briefly review the projection-based model reduction methods with nonlinear manifold augmentation. Then in Section~\ref{sec:hybrid}, we present the hybrid ROM. Section~\ref{sec:non-intrusive} proposes a non-intrusive variant of the hybrid ROM. Section~\ref{sec:results} shows the numerical results. Conclusions are drawn in the end. 

\section{Background }
\label{sec:POD-greedy}

We consider surrogate modeling of parametric time-dependent problems in the general nonlinear form as below,
\begin{equation}
\begin{array}{rl}
	\label{eq:fom}
E(\mu) \frac{d}{dt}{\bf u}(t, \mu) &= A(\mu) {\bf u}(t, \mu) + {\bf f}({\bf u}(t, \mu)) + B(t, \mu)  \qquad {\bf u}(0, \mu) = {\bf u}_{0}(\mu)
\end{array}
\end{equation}
where $t\in [0, T]$ and  $\mu \in \mathcal{P} \subset \mathbb R^{p}$, $\mathcal{P}$ is the parameter domain. ${\bf u}(\mu)\in \mathbb R^n$ and ${\bf E}(\mu), {\bf A}(\mu) \in \mathbb R^{n\times n}$ are the linear operators depending on the parameter $\mu$. $B: (t, \mu) \mapsto \mathbb R^{n}$ is the source term, and  ${\bf f}: \mathbb R^N \times \mathcal P \mapsto \mathbb R^{N}$ is the nonlinear operator. 
The system in~(\ref{eq:fom}) is usually the result of numerical discretization of partial differential equations (PDEs), which is named as the full-order model (FOM). The number of equations $n$ in~(\ref{eq:fom}) is often very large to ensure high-resolution of the physical dynamics. 

\subsection{The reduced basis method}
\label{sec:RBM}
The RBM finds a ROM of the FOM in~(\ref{eq:fom}) with small number $r\ll n$ of equations using a POD-greedy process. We briefly describe the process in Algorithm~\ref{alg:POD-greedy}. 
  \begin{algorithm}[H]
  \caption{POD-greedy algorithm for parametric time-dependent systems}
  \begin{algorithmic}[1]
  \REQUIRE  \hspace{0.15cm}    $\mathcal P_{train},tol_{RB}(<1)$, $r_{\max}$
  \ENSURE   Projection matrix: $\tilde V=[V_1,\ldots,V_{\tilde r}] \in R^{n\times \tilde r}$
  \STATE Initialization: $s=0$, $\tilde V=[\,]$,
  $\mu^{\star}=\mu^1$, $\eta(\mu^{\star})=1$.
  \WHILE{  the error $\eta(\mu^{\star})>tol_{RB}$ or $\tilde r\leq r_{\max}$} 
  \STATE Compute the snapshot matrix $X=[{\bf u}(t_1, \mu^{\star}), \ldots, {\bf u}( t_{N_t}, \mu^{\star})]$.
    \STATE   $\mathcal W=\textrm{colspan}\{V\}$. 
\STATE
IF $\tilde V\neq$ [\ ] do

$X=X-\tilde V \tilde V^TX$. (Remove the $\tilde V$-dependent part from $X$ and update $X$.)
 
ENDIF
\STATE Do singular value decomposition (SVD): $X=U_x\Sigma_x V_x^T$, 
$V_{s+1}:=U_x(:,1:r_s)$. 
\STATE Update $\tilde V$: $\tilde V=[\tilde V,V_{s+1}]$,
     \STATE  $s=s+1$.
  \STATE  Find $\mu^{\star}:= \mbox{arg}\max\limits_{\mu \in \mathcal P_{train}}\eta(\mu)$.
  \ENDWHILE
  \end{algorithmic}
  \label{alg:POD-greedy}
  \end{algorithm}
In Algorithm~\ref{alg:POD-greedy},  the truncation order $r_s$ is determined by $$\left(\sum\limits_{j=r_s+1}^d \sigma_j\right)/ \left(\sum\limits_{j=1}^d \sigma_j\right)<tol_{svd},$$ 
where $\sigma_j, j=1\ldots, d$ are the singular values on the diagonal of $\Sigma$, and $tol_{svd}$ is a small value, e.g., 0.01. When $r_s=1$, it is the standard POD-greedy algorithm, where only the first dominant singular vector is taken to enrich $\tilde V$, i.e., $V_{s+1}:=U(:,1)$. For convection-dominated problems, the same parameter sample might be repeatedly chosen to catch the resolution of the corresponding solution. Choosing multiple dominant left singular vector at each iteration might reduce such repetitions to a large extend. The tolerance $tol_{svd}$ is used to determine the number of dominant singular vectors. We use the 2-norm of the residual vector obtained from simulating the ROM as the error estimator $\eta(\mu)$, i.e., $\eta(\mu)=\|{\bf r} (T, \mu)\|_2$, the magnitude of the residual vectors at the final time point. 
Here, ${\bf r} (t, \mu):=E(\mu) d \tilde {\bf u}(t, \mu)/dt-\tilde {\bf u}(t, \mu) -A(\mu) \tilde {\bf u}(t, \mu)- {\bf f} (\tilde {\bf u}(t, \mu))-B(t, \mu), \forall t \in [0, T]$, and $\tilde {\bf u}(t, \mu)$ is the approximate solution. 
A more accurate error estimator~\cite{morFenCB24} can also be used here. Since this is not the focus of this work, we postpone it to a future work. 

\subsection{Galerkin-based ROM}

After the projection matrix $\tilde V$ is computed, the RBM usually constructs a ROM based on Galerkin projection,
\begin{equation}
\label{eq:rom}
\begin{array}{rcl}
\tilde V^TE(\mu)\tilde V d{\bf z}(t, \mu)/dt&=& \tilde V^TA(\mu) \tilde V {\bf z}(t, \mu)+ \tilde V^T{\bf f} (\tilde V {\bf z}(t, \mu))+\tilde V^TB(t, \mu).
\end{array}
\end{equation}
The complexity of the ROM in~(\ref{eq:rom}) needs to be further reduced by hyper-reduction methods, e.g. the empirical interpolation method (EIM)~\cite{morBarMNetal04}, the discrete empirical interpolation method (DEIM)~\cite{morChaS10}, etc. We employ the EIM in this work. Finally, we can get the following ROM with EIM-reduced complexity,  
\begin{equation}
\label{eq:rom_EI}
\begin{array}{rcl}
\tilde V^TE(\mu)\tilde V d{\bf z}(t, \mu)/dt&=& \tilde V^TA(\mu)\tilde V {\bf z}(t, \mu)+ \tilde V^T U_{EI} \beta(t, \mu)+\tilde V^TB(t, \mu),
\end{array}
\end{equation}
where ${\bf f} (\tilde V {\bf z}(t, \mu))$ is replaced by the linear expansion of the EI basis vectors, i.e. the column vectors in $U_{EI}$. The coefficient vector $\beta(t, \mu)=(P^TU_{EI})^{-1}P^T{\bf f}(\tilde Vz(t, \mu))$. $P^T{\bf f}(\tilde Vz(t, \mu))$ is equivalent to evaluating the nonlinear vector ${\bf f}(\tilde Vz(t, \mu))$ only at those entries indicated by $P^T$, so that computing ${\bf f} (\tilde V {\bf z}(t, \mu ))$ can be independent of the large dimension $n$. The final approximate solution $\tilde {\bf u}(t, \mu)$ of the FOM can be computed as $\tilde {\bf u}(t, \mu)=\tilde V {\bf z}(t, \mu)$.  We call the ROM in (\ref{eq:rom_EI}) the Galerkin-ROM (G-ROM).  In Algorithm~\ref{alg:POD-greedy}, we need to compute the error estimator $\eta(\mu)$ for obtaining the final G-ROM.  In fact, $\eta(\mu)$ is computed from the approximate solution produced by the EIM-reduced G-ROM that is updated at each iteration whenever the matrix $\tilde V$ is updated.

\subsection{Least-squares Petrov-Galerkin-based ROM }
Instead of the G-ROM, one may use Petrov-Galerkin ROM to improve stability~\cite{morGriFY20} or accuracy~\cite{morCarBF11}. Here, we review the least-squares-Petrov-Galerkin-based ROM (LSPG-ROM) proposed in~\cite{morCarBF11}. 
With the matrix $\tilde V$, one gets the residual induced by the approximate solution $\tilde {\bf u}(t, \mu)=\tilde V{\bf z}(t, \mu)$ as
\begin{equation}
\label{eq:resi}
\begin{array}{rcl}
{\bf r}(\tilde V{\bf z}(t, \mu))=E(\mu) d \tilde V {\bf z}(t, \mu)/dt- A(\mu) \tilde V{\bf z}(t, \mu)-{\bf f}(\tilde V{\bf z}(t, \mu))-B(t, \mu)
\end{array}
\end{equation}
When time discretization is applied to (\ref{eq:resi}), we get the residual ${\bf r}^i (\tilde{\bf u}(t_i, \mu))$ at each $t_i$. Expression of ${\bf r}^i$ depends on the selected time discretization scheme. Using backward finite difference discretization in the time domain, we have
\begin{equation}
\label{eq:resi1}
\begin{array}{rcl}
{\bf r}^i (\tilde V{\bf z}(t_i, \mu))=\frac{1}{\Delta t}E(\mu)\tilde V{\bf z}(t_i, \mu)-\frac{1}{\Delta t}E(\mu)\tilde V{\bf z}(t_{i-1}, \mu)- A(\mu) \tilde V{\bf z}( t_i, \mu)-{\bf f}(\tilde V{\bf z}(t_i, \mu))-B(t_i, \mu),
\end{array}
\end{equation}
where $\Delta t$ is the time step size, assuming that a uniform step size is used.  
Note that ${\bf r}^i \in \mathbb R^n$ is nonzero and is still in the original large-dimensional space.  The reduced state vector ${\bf z}(t_i, \mu)$ can be computed by solving the following minimization problem, 
$${\bf z}(t_i, \mu)=\min\limits_{\tilde {\bf z}\in \mathbb R^{\tilde r}}\|{\bf r}^i (\tilde V\tilde {\bf z}(t_i, \mu))\|_2$$
It is proved in~\cite{morGriFY20} that the above minimization problem is equivalent to solving the following LSPG-ROM,
\begin{equation}
\label{eq:LSPG-ROM}
W({\bf z}(t_i, \mu))^T {\bf r}^i (\tilde V{\bf z}(t_i, \mu))=0, 
\end{equation}
where the left projection matrix $W$ spanning the testing space, is computed as
$$W({\bf z}(t_i, \mu))=J(\tilde{\bf u}(t_i, \mu)) \tilde V.$$
The Jacobian matrix is the Jacobian of ${\bf r}^i$ w.r.t the approximation solution $\tilde {\bf u}(t_i, \mu)$, i.e.,
 $$J(\tilde{\bf u}(t_i, \mu))=\frac{\partial {\bf r}^i}{\partial \tilde {\bf u}} (\tilde {\bf u}(t_i, \mu)).$$
The LSPG-ROM~(\ref{eq:LSPG-ROM}) is now in the low-dimensional space $\mathbb R^{\tilde r}$ and can be solved by the Newton's method at each $t_i$. At each Newton iteration $j$, the Jacobian depends on the Newton approximation solution $\tilde {\bf u}^j(t_i, \mu) \approx \tilde {\bf u}(t_i, \mu)$ and can be computed as,
$$J^{i,j} =\frac{\partial {\bf r}^i}{\partial \tilde {\bf u}} (\tilde {\bf u}^j(t_i, \mu))=\frac{\partial {\bf r}^i}{\partial \tilde {\bf u}} (\tilde V{\bf z}^j (t_i, \mu)).$$
More detailed analysis of the method can be found in~\cite{morCarBF11,morGriFY20}. It is clear that the complexity of solving the LSPG-ROM needs to be further reduced by hyper-reduction, as evaluations of ${\bf r} (\tilde V \tilde {\bf z}^j(t_i, \mu))$ and its Jacobian matrix at each Newton iteration and each $t_i$ still have complexity depending on the large dimension $n$. Similarly, we apply the EIM method to (\ref{eq:LSPG-ROM}), in particular to both ${\bf r}(\tilde V \tilde {\bf z}^j(t_i, \mu))$ and $J^{i,j}$. Finally, the LSPG-ROM introduced much more computational complexity in both the offline computation and the online simulation as compared to the G-ROM. 

\subsection{ROMs with nonlinear manifold augementation}
\label{sec:NManifold}
The G-ROM or the LSPG-ROM is  of small size $\tilde r$ when the solution manifold can be well presented by a linear subspace. This is nevertheless not true for convection-dominated problems, especially for nonlinear hyperbolic conservation laws. On the contrary, we may need a ROM with relatively large $r$ for sufficient accuracy. In order to keep the size $\tilde r$ small, nonlinear manifold methods are proposed in~\cite{morBarFM2023,morParTHetal26,morBarC22,morAghBMetal26}, where instead of using only $\tilde V$ as the trial space for approximating ${\bf u}(t, \mu)$, a nonlinear manifold is used to approximate the subspace spanned by $V_2$. Here, $\tilde V=(V_1, V_2)$ is divided into two matrices with $V_1 \in \mathbb R^{n\times r}, r\ll \tilde r$, so that  $\tilde {\bf u}(t, \mu)=\tilde V {\bf z}(t, \mu)$ in~(\ref{eq:rom_EI}) or in~(\ref{eq:LSPG-ROM}) can be written as $\tilde {\bf u}(t, \mu)=V_1 {\bf z}_1(t, \mu)+V_2 {\bf z}_2(t, \mu)$. 

The methods of nonlinear manifold approximates ${\bf z}_2(t, \mu)$ with a nonlinear function of ${\bf z}_1(t, \mu)$, so that the solution is approximated as ${\bf u}(t, \mu)\approx V_1 {\bf z}_1(t, \mu)+ V_2 g({\bf z}_1(t, \mu))$. The G-ROM or the LSPG-ROM can be obtained by replacing $\tilde {\bf u}(t, \mu)=\tilde V {\bf z}(t, \mu)$ with $\bar {\bf u}(t, \mu)=V_1 {\bf z}_1(t, \mu)+ V_2 g({\bf z}_1(t, \mu))$. The nonlinear operator $g({\bf z}_1(t, \mu))$ can be a quadratic function of ${\bf z}_1(t, \mu)$~\cite{morBarC22,morAghBMetal26}, a function learned by radial basis interpolation~\cite{morParTHetal26}, or a function represented by a neural network model~\cite{morBarFM2023}. Then a ROM with the unknown solution vector ${\bf z}_1(t, \mu) \in \mathbb R^r, r\ll \tilde r$ is solved to get the approximate solution $\bar {\bf u}(t, \mu)$. It can be expected that the resulting ROM, especially the LSPG-ROM, requires involved offline and online computations when they are further hyper-reduced by EIM or other hyper-reduction methods. Moreover, since $g({\bf z}_1(t, \mu))$ is a function that approximates ${\bf z}_2$, the $r$-dimensional ROMs with nonlinear manifold augmentation $\bar {\bf u}(t, \mu)=V_1 {\bf z}_1(t, \mu)+ V_2 g({\bf z}_1(t, \mu))$ cannot be more accurate than the $\tilde r$-dimensional ROMs with the linear subspace approximation $\tilde {\bf u}(t, \mu)=\tilde V {\bf z}(t, \mu)$. 
In the next section, we propose a hybrid ROM in which the model component takes the form of a low-dimensional G-ROM (LSPG-ROM) and is kept as simple and standard as possible. The error of the model component is then captured by physics-informed deep learning. The hybrid ROM breaks the barrier of linear-subspace-limited-accuracy of the ROMs with nonlinear manifold augmentation. 
\section{The proposed hybrid ROM}
\label{sec:hybrid}
\begin{figure}[!h]
\centering
\includegraphics[width=100mm]{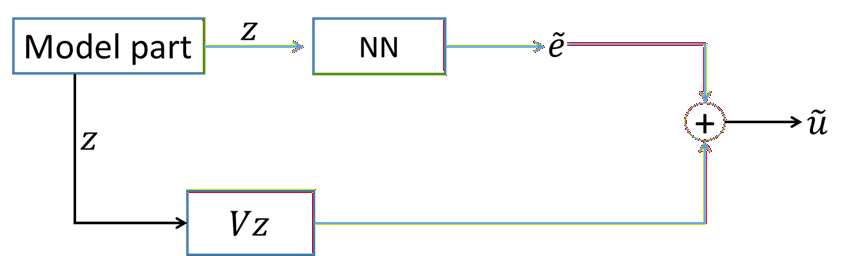}
\caption{Structure of the hybrid ROM.}
\label{fig:hybridROM}
\end{figure}

As already mentioned in the Introduction, the hybrid ROM is of the form,
\begin{equation}
\label{eq:ROM-hybrid}
\tilde {\bf u}(t,\mu)=V {\bf z}(t,\mu)+\tilde {\bf e}({\bf z}(t, \mu)),
\end{equation}
where ${\bf z}(t, \mu) \in \mathbb R^r$ is computed from the ROM in~(\ref{eq:rom_EI}), with a small value $r$. Figure~\ref{fig:hybridROM} shows the structure of the hybrid ROM. 
The matrix $V \in \mathbb R^{n\times r}$ can be computed in the following different ways:
\begin{itemize}
\item It is directly computed using the POD-greedy algorithm by requiring $r_{\max}=r$. 

\item It can also be computed from the POD-greedy algorithm, with $r_{\max} \gg r $. After POD-greedy converged, the snapshot matrices corresponding to the parameters selected by POD-greedy are used to constitute a matrix $X_s$, then $V$ is obtained as the first $r$ left singular vectors of the SVD of $X_s$, i.e., $V=U_s(:, 1:r)$, $[U_s, \sim]$=SVD( $X_s, $ `econ'). Here, we use the MATLAB notation. In this way, $V$ could capture richer information in the parameter space.  

\item It is computed from the SVD of $X_{all}$. $X_{all}$ is the matrix composed of the snapshots corresponding to all the training parameter samples in $\mathcal P_{train}$. For very large problems, this might occupy too much memory when the number of training samples is not small. 
\end{itemize}
The model component $V {\bf z}(t,\mu)$ of the hybrid ROM in~(\ref{eq:ROM-hybrid}) corresponds to the a standard G-ROM or LSPG-ROM. I.e., the reduced states ${\bf z}(t,\mu)$ are the solution vector of a standard G-ROM or LSPG-ROM. The other component of the hybrid ROM, $\tilde {\bf e}({\bf z}(t, \mu))$:  $\mathbb R^r \mapsto \mathbb R^n$, is a function of ${\bf z}(t, \mu)$. It is learned from the NN,
$$\tilde {\bf e}({\bf z}(t, \mu))=f_{NN} ({\bf z}(t, \mu)).$$ 
$\tilde {\bf e}({\bf z}(t, \mu))$ is an approximation of the error ${\bf u}(t, \mu)-V {\bf z}(t,\mu)$ contributed by the model component of the hybrid model. 
The input of the NN is the reduced (latent) states ${\bf z}(t, \mu)$ computed from the model component, rather than the time-parameter pair $(t, \mu)$. This brings physical information into the NN. Finally, $f_{NN} ({\bf z}(t, \mu))$ is the NN component of the hybrid ROM and the final approximate solution $\tilde{\bf u}(t, \mu)$ is the sum of the model component approximation $V{\bf z}(t, \mu)$ and the error $\tilde {\bf e}({\bf z}(t, \mu))$ learned from NN.
\subsection{The NN structure}
\label{sec:NNstructure}
\begin{figure}[!h]
\centering
\includegraphics[width=100mm]{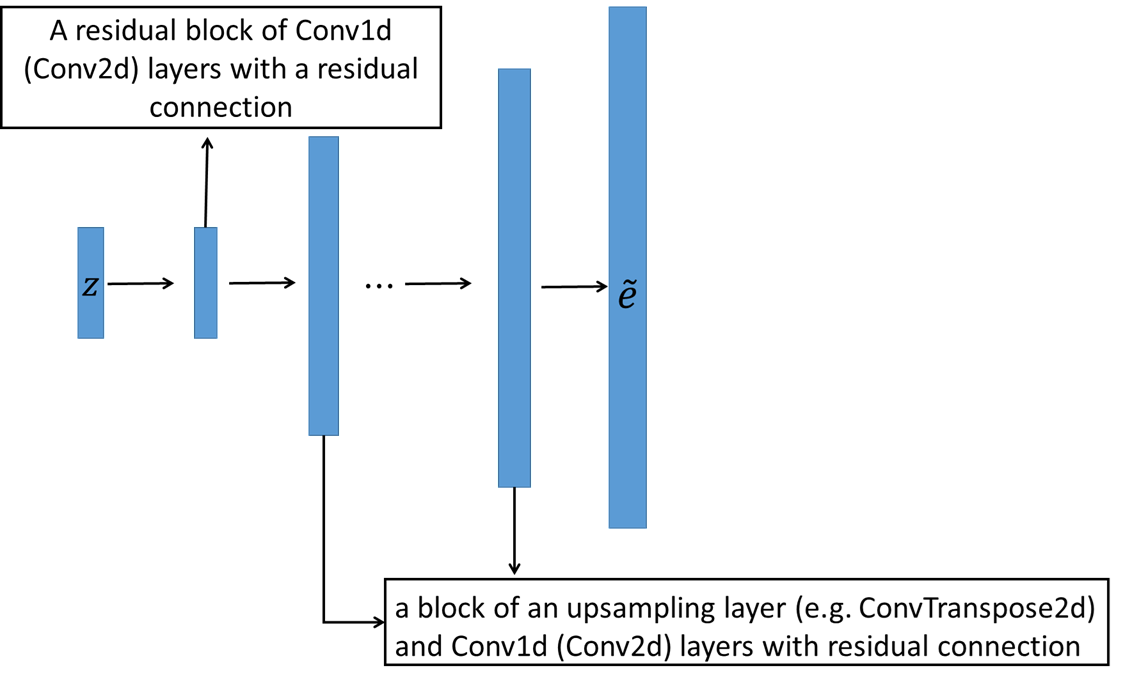}
\caption{The NN structure of the hybrid ROM}
\label{fig:1dNN}
\end{figure}

Depending on the problem considered, the structure of the NN, $f_{NN}$, can be different. In this work, we propose to use a convolutional neural network (CNN) with residual connections between the hidden layers. The convolutional neural network is mainly used to learn the spatial relation of the data. The residual connection is important for avoiding vanishing gradient and improve representation quality. In addition to Conv1d (Con2d) layers with residual connections, we use upsampling layers to decode the short input vector into the long error vector in the physical domain. More specifically, 
\begin{itemize}
\item For 1D problems, $f_{NN}$ includes one residual block and multiple upsampling blocks, each containing  an embedded residual block. All these constitute a decoder network mapping the reduced state vector ${\bf z}(t, \mu)$ to the error vector $\tilde {\bf e}(t, \mu)$ in the original physical domain. The general structure is illustrated in Fig.~\ref{fig:1dNN}. Between the input and the output layers, there are several upsampling blocks composed of upsampling layers, Conv1d layers and skip connections. The number of the Conv1d layers in each block can be different for different problems. We name this CNN with residual connection as 1dCNNResi.

\item For 2D problems, $f_{NN}$ is a CNN containing a few upsampling (ConTranspose2d) layers, each being concatenated with two sequential residual blocks. Each residual block, in turn contains a few Conv2d layers and a skip connection. It is a decoder network mapping the reduced state vector ${\bf z}(t, \mu)$ to the error vector $\tilde {\bf e}(t, \mu)$ in the original physical domain. If there are two PDEs in the original problems, then the error vector is separated  into two channels, each representing the solution error of a single scalar PDE. The general structure is similar as that in Fig.~\ref{fig:1dNN}. We use ConvTranspose2d as the upsampling layer for the 2D problems.  This $f_{NN}$ for 2D problems is named as 2dCNNResi. 

\item Very large FOMs often arise from 2D problems, then the 2dCNNResi has many NN parameters, making the training time relatively long. The original data can be pre-processed by POD compression, so that only the POD coefficient vectors with fare less entries than the original error vectors are learned by the 2dCNNResi.  This $f_{NN}$ is named as POD-2dCNNResi. More specifically, 
\begin{equation}
\label{eq:PODcoefficients}
\tilde {\bf e}({\bf z}(t, \mu))=V_0 \hat {\bf e}({\bf z}(t, \mu)),
\end{equation}
where  $\hat {\bf e}({\bf z}(t, \mu))\in \mathbb R^{r_0}, r\leq r_0 \ll n$ is learned by a 2dCNNResi. The output of this 2dCNNResi is a vector with length $r_0$ that is much shorter than the length-$n$ of output produced by the 2DCNNResi for the previous case. The matrix $V_0$ includes the dominant left singular vectors of the error snapshot matrix $X_e=[{\bf e}(t_0, \mu_1), \ldots, {\bf e}(t_{n_t-1}, \mu_1) ), \ldots, {\bf e}(t_0, \mu_{n_\mu}), {\bf e}(t_{n_t-1}, \mu_{n_\mu})]$ at the training parameter samples and training time instances. After pre-processing, the 2dCNNResi will contain much reduced number of NN parameters, so that the training time can be largely shortened. 
More detailed explanations for 1dCNNResi and 2dCNNResi can be found in Section~\ref{sec:results} for each example there. Additional figures illustrating the details of the corresponding NN structure are presented in~\ref{sec:appendix}.

\end{itemize}


%
%
%
%
%
%
%
%
%

\section{A non-intrusive variant}
\label{sec:non-intrusive}
Based on the NN structure of the hybrid ROM, we propose a non-intrusive variant: POD-FFNN-e-decoder. Instead of computing the reduced state vectors ${\bf z}(t, \mu)$ from the projection based ROM, we learn them from a feedforward neural network (FFNN). The error vector is still learned from the NNs presented in Section~\ref{sec:hybrid}, but with inputs being the outputs of the FFNN. That means, the reduced states $\tilde {\bf z}(t, \mu)$ learned from FFNN are fed into the NNs in Section~\ref{sec:hybrid} to reconstruct the error vectors $\tilde {\bf e}(t, \mu)$, so that each NN presented in Section~\ref{sec:hybrid} acts as an error decoder (e-decoder). 

\begin{figure}[!h]
\centering
\includegraphics[width=100mm]{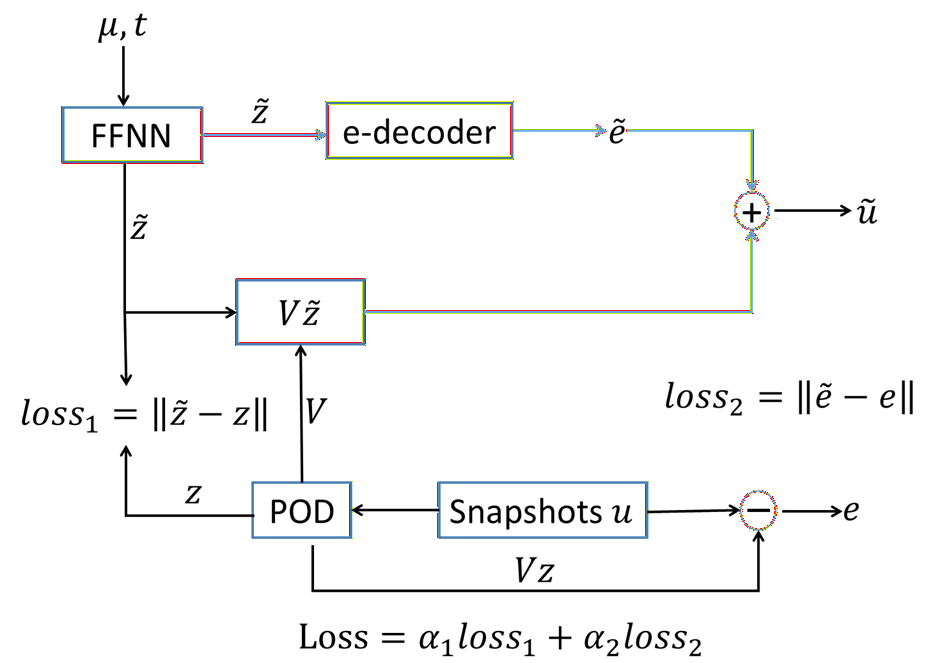}
\caption{Structure of the POD-FFNN-e-decoder}
\label{fig:non-intrusive}
\end{figure}

Figure~\ref{fig:non-intrusive} presents the flowchart of proposed POD-FFNN-e-decoder. 
Let the $X$ be the matrix composed of the solution snapshots, i.e., $$X=[{\bf u}(t_0, \mu_1), \ldots, {\bf u}(t_{n_{t-1}}, \mu_1), \ldots, {\bf u}(t_0, \mu_{n_\mu}), \ldots, {\bf u}(t_{n_{t-1}}, \mu_{n_\mu})],$$ 
 where $\mu_1, \ldots, \mu_{n_\mu} \in \mathcal P_{train}$. We apply SVD to $X$ and get the matrix $U_x \in \mathbb R^{n\times n_t n_\mu}$ of the left singular vectors. The first $r$ dominant singular vectors are used to compute the POD-reduced states, 
\begin{equation}
\label{eq:solution-PODcoefficients}
\begin{array}{rcl}
Z&=&V^T X, V=U_x(:, 1:r)\\
Z&:=&[{\bf z}(t_0, \mu_1), \ldots, {\bf z}(t_{n_t}, \mu_1), \ldots, {\bf z}(t_0, \mu_{n_\mu}), \ldots, {\bf z}(t_{n_{t-1}}, \mu_{n_\mu})]
\end{array}
\end{equation} 
in Fig.~\ref{fig:non-intrusive}, which are the training output data for the FFNN and training input data for the e-decoder, 1dCNNResi, or 2dCNNresi. The training output data for the e-decoder are
computed as 
\begin{equation}
\label{eq:errorsnapshots_edecoder}
\begin{array}{rcl}
X_e&=&X-VZ, \\ 
X_e&:=&[{\bf e}(t_0, \mu_1), \ldots, {\bf e}(t_{n_t}, \mu_1), \ldots, {\bf e}(t_0, \mu_{n_\mu}), \ldots, {\bf e}(t_{n_{t-1}}, \mu_{n_\mu})].
\end{array}
\end{equation} 

The FFNN is trained together with the decoder. Moreover, we have two loss functions: one defines the reconstruction accuracy of the decoder for the error vectors, the other measures the learning accuracy of FFNN. The $loss_1$ in Fig.~\ref{fig:non-intrusive} for the FFNN is evaluated via the mean squared error (MSE) between  $\tilde {\bf z}(t, \mu)$ and the reduced states ${\bf z}(t, \mu)$ in~(\ref{eq:solution-PODcoefficients}). The $loss_2$ in Fig.~\ref{fig:non-intrusive} for the e-decoder is the MSE between ${\bf e}(t, \mu)$ in~(\ref{eq:errorsnapshots_edecoder}) and the predicted value $\tilde {\bf e}(t, \mu)$. The final loss is a weighted sum of $loss_1$ and $loss_2$, where $\alpha_1$ and $\alpha_2$ control the contribution of each loss to the final loss. 
\begin{remark}
Note that the proposed hybrid ROM and its non-intrusive variant POD-FFNN-e-decoder share a similar framework as the one in~\cite{morCodIDetal25}: $\hat u({\bf x},\mu)=V \hat z({\bf x},\mu)+\hat z({\bf x}, \mu)^T C({\bf \Phi}, \mu) \hat z({\bf x}, \mu)$, where the error correction term is assumed to be a quadratic function of $\hat z({\bf x}, \mu)$, and $\hat z({\bf x}, \mu)$ is obtained from a non-intrusive method based on POD and radial basis function interpolation. The columns of ${\bf \Phi}\in \mathbb R^{n\times r}$ are the $r$ dominant POD basis vectors.  A neural operator is constructed in~\cite{morCodIDetal25} to learn the tensor $C({\bf \Phi}, \mu)$. However, the proposed ROM $\hat u({\bf x},\mu)$ is not time-dependent. Whether the method is applicable to transient solution of time-dependent problems, is unclear.  
\end{remark}

\section{Numerical results}
\label{sec:results}
This section presents the numerical tests on a 1D inviscid Burgers' equation with a single parameter on the initial condition, and on 2D inviscid Burgers' equations with two parameters on the boundary and in the source term, respectively. 
The model components of the hybrid ROMs for the two problems as well as the standard RBM method are implemented on a {\it single} core of a compute-server with two AMD EPYC 7763 64-core processors and 1TB main memory. The NN component of the hybrid ROM and POD-FFNN-e-decoder are trained and tested on a Ubuntu desktop with 31GiB, 12th Gen Intel(R) Core(TM) i5-12600K x 16 processor, and GPU. 

We use the relative errors defined as below to measure the accuracy of the proposed hybrid ROM. 
$$\epsilon(\mu)=\frac{\sum\limits_{i=1}^{N_t}\|{\bf u}(t_i, \mu)-\tilde{\bf u}(t_i, \mu)\|_2}{\sum\limits_{i=1}^{N_t}\|{\bf u}(t_i, \mu)\|_2}.$$
$$\epsilon_{\max}=\max\limits_{\mu \in {\mathcal P}_{test}} \epsilon(\mu),$$
where $\mathcal P_{test}$ is the set of testing parameter samples.
We also present the wall clock time for offline training and online prediction. In the relevant tables below, h means hour, m means minutes, and s means seconds.

\subsection{1D inviscid Burgers' equation}
The 1D inviscid Burgers' equation is given by:
\begin{equation}
\frac{\partial u}{\partial t} + u \frac{\partial u}{\partial x} = 0, \quad x \in [0,L], \; \quad t \in [0,T],
\end{equation}
with a parametrized initial condition:
\begin{equation}
u(x,0;\mu) =
\begin{cases} 
1 + \frac{\mu}{2} \sin(2 \pi x - \pi/2) + 1, & 0 \leq x \leq 1, \mu\in[0.9, 1.1] \\
1, & \text{otherwise},   
\end{cases}
\end{equation}
and periodic boundary condition: $$u(0,t, \mu)=u(L,t,\mu).$$ 
We follow the finite difference method in~\cite{morKimCWZ22} to obtain the spatially discretized FOM~(\ref{eq:fom}) with $L=2, T=0.5, n=1000$. The parameter domain is also the same interval as in~\cite{morKimCWZ22}. The 1D space is discretized into 1000
equal elements, resulting in 1001 grid points.  In particular, the backward finite difference scheme is used: $\frac{\partial u}{\partial x}=\frac{u_i-u_{i-1}}{\Delta x}, \Delta x=L/n$. With the periodic boundary condition, we finally have a state vector ${\bf u}$ of $n=1000$ entries in the FOM, corresponding to the values of the solution $u(x, t, \mu)$ at 1000 spatial grid points. In this section, we compare four surrogate models, the RBM G-ROM, the hybrid ROM, POD-FFNN-e-decoder and an existing CAE-FFNN model. All surrogates are trained using the same parameter and time samples. 

\subsubsection{Training the model component of the hybrid ROM}
The model component of the hybrid ROM is obtained by the POD-greedy algorithm with the following settings:
\begin{itemize}
    \item Tolerance for the greedy algorithm: $tol_{RB} = 0.05$  
    \item Tolerance for the EIM: $tol_{EIM}=10^{-5}$.   
   \item $r_s=1$ in Step 6. 
    \item Training parameter sample set: $\mathcal P_{train} =\{\mu_1, \ldots, \mu_{n_\mu}\}=\{0.9, 0.95, 1, 1.05, 1.1\}$  
\end{itemize}

At each iteration of Algorithm~\ref{alg:POD-greedy}, the error estimator is computed from the G-ROM~(\ref{eq:rom_EI}). Finally, the POD-greedy algorithm gives a projection matrix $\tilde V$ with $\tilde r=44$ columns, by choosing only the snapshots at $\mu=0.9, 1.1$. Matrix $\tilde V$ is used to construct a RBM G-ROM with reduced size $\tilde r=44$ for comparison. We then take the first $r=5$ columns from $\tilde V$ to form the matrix $V$ for the model component of the hybrid ROM, $V=\tilde V(:,1:5)$. Applying this $V$ to~(\ref{eq:rom_EI}) yields the model component of the hybrid ROM. By simulating the model component, we can get the reduced state vector ${\bf z}(t,\mu) \in \mathbb R^r$ at any time and $\mu$.

\subsubsection{Training the NN component of the hybrid ROM}
The NN component, 1dCNNResi, for the 1D Burgers' equation includes one residual block depicted in Fig.~\ref{fig:resi_block}, where the output of the residual block can be written as $\mathrm{ResiBlock}({\bf x})=\mathrm{ReLU}({\bf x}+$ conv1d(ReLU(conv1d (${\bf x}$))). There is another residual block that is embedded in the upsampling block, where the input $x$ is fed into a Conv1d layer, and then skip connection is added to the final Conv1d layer. The final output of this residual block is named as $\mathrm{ProjResiBlock}({\bf x})$ and can be written as $\mathrm{ProjResiBlock}({\bf x})$=$\mathrm{ReLU}$(Con1d(${\bf x}$)+ conv1d(ReLU(conv1d (${\bf x}$))). It is illustrated in Fig.~\ref{fig:Projresi_block}. The Upsampling block in the CNN decoder contains a Upsampling layer and the ProjResiBlock, as shown in Fig.~\ref{fig:up_block}. Finally, the structure of the 1dCNNResi is given in Fig.~\ref{fig:CNNresi}, which is a sequential connection of the residual block ResiBlock with 6 upsampling blocks. 

To get the error snapshots for training the 1dCNNResi, we simulate the FOM at the training samples of $\mu \in  \mathcal P_{train}$ using backward Euler method for time integration at $n_t=500$ time steps.  At each $t_i, i=1, \ldots, n_t-1$, the Newton's method is used to get the solution ${\bf u}(t_i, \mu)$.  The same time-integration scheme is applied to solve the model component and get the approximate solution $\hat {\bf u}(t_i, \mu)=V{\bf z}(t_i, \mu)$. Once both ${\bf u}(t_i, \mu)$ and $\hat {\bf u}(t_i, \mu)$ are obtained for each training sample of $\mu$, we can get the error snapshots ${\bf e}(t_i, \mu_j), i=0, \ldots, n_t-1, \mu_j \in \mathcal P_{train}$ for training 1dCNNResi. The complete training dataset for the 1dCNNResi is $\{{\bf z}(t_i, \mu_j)$, ${\bf e}(t_i, \mu_j)\}, i=0, \ldots, n_t-1, \mu_j \in \mathcal P_{train}\}$. Afterwards, 1dCNNResi is trained for 1000 epochs.

\begin{remark}
Since the FOM needs to be simulated at all $\mu \in  \mathcal P_{train}$ to get the snapshots of ${\bf f}({\bf z}(t, \mu))$ for implementing hyper-reduction, the solution snapshots ${\bf u}(t_i, \mu_j), i=0, \ldots, n_t-1, \mu_j \in \mathcal P_{train}$ are obtained as a byproduct of the EIM method, and are used for free when implementing the POD-Greedy algorithm as well as training the neural network. That means, we don't need to solve the FOM at the selected samples of $\mu$ at Step 3 of Algorithm~\ref{alg:POD-greedy}. The same solution snapshots are then also used for generating training data of the neural networks 1dCNNResi, POD-FFNN-e-decoder and CAE-FFNN. 
\end{remark}

\subsubsection{Training the POD-FFNN-e-decoder}
The POD-FFNN-e-decoder for the 1D Burgers' equation is composed of a FFNN and 1dCNNResi as the e-decoder. These two networks are trained simultaneously with the loss defined in Fig.~\ref{fig:non-intrusive}. The training data for FFNN are computed following~(\ref{eq:solution-PODcoefficients}) with $r=5$. The training data for 1dCNNResi are the error snapshots computed from solution snapshots and the POD reconstructed approximate snapshots, which are computed as in~(\ref{eq:errorsnapshots_edecoder}). The latent space dimension $r=5$ coincides with the reduced state dimension of the hybrid-ROM. The POD-FFNN-e-decoder is also trained for 1000 epochs.

We also compare the hybrid ROM and the  POD-FFNN-e-decoder with an existing non-intrusive MOR method CAE-FFNN in~\cite{morFreetal21}. During the offline stage, a CAE and a feedforward neural network (FFNN) are trained with the {\it solution} snapshots that are also used to generate the training data for 1dResiCNN. During the online stage, the decoder of the CAE plus the FFNN are used to predict the solution at any given time instance and given parameter sample. In order to have a fair comparison, the decoder of the CAE has the same structure of our 1dCNNResi, and the encoder has the corresponding downsampling structure with the same skip connections. The same FFNN used in POD-FFNN-e-decoder is also used in CAE-FFNN. To be in agreement with 1dCNNResi and POD-FFNN-e-decoder, the CAE compresses the solution data into a latent space of dimension $r=5$. The CAE-FFNN is trained for 1000 epochs. 

\subsubsection{Prediction results}
\begin{figure}[!h]
\centering
\includegraphics[width=100mm]{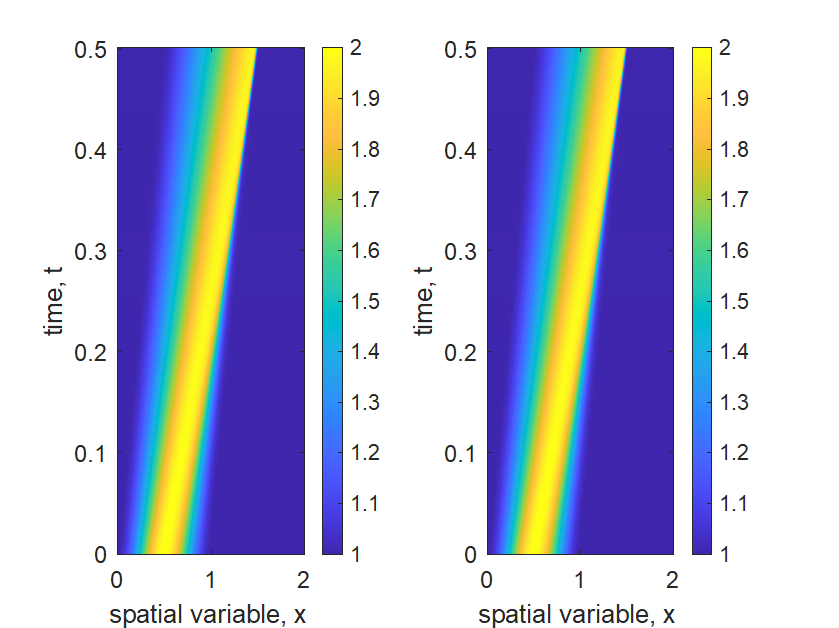}
\caption{The FOM solution at $\mu^*=1.08$ (left) and the hybrid ROM at $\mu^*$ (right).}
\label{fig:solution_comp1}
\end{figure}

\begin{figure}[!h]
\centering
\includegraphics[width=80mm]{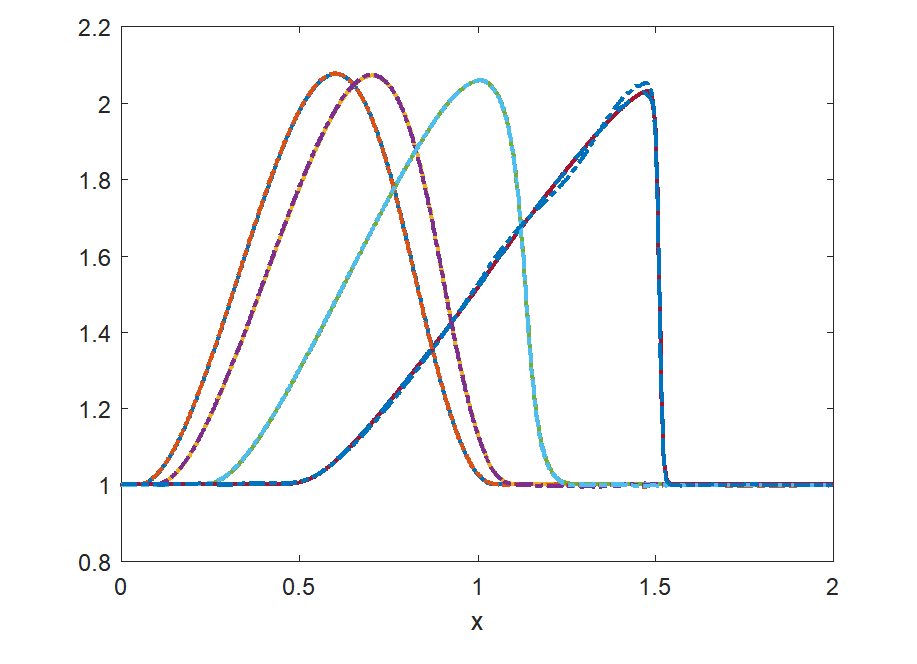} \quad 
\includegraphics[width=80mm]{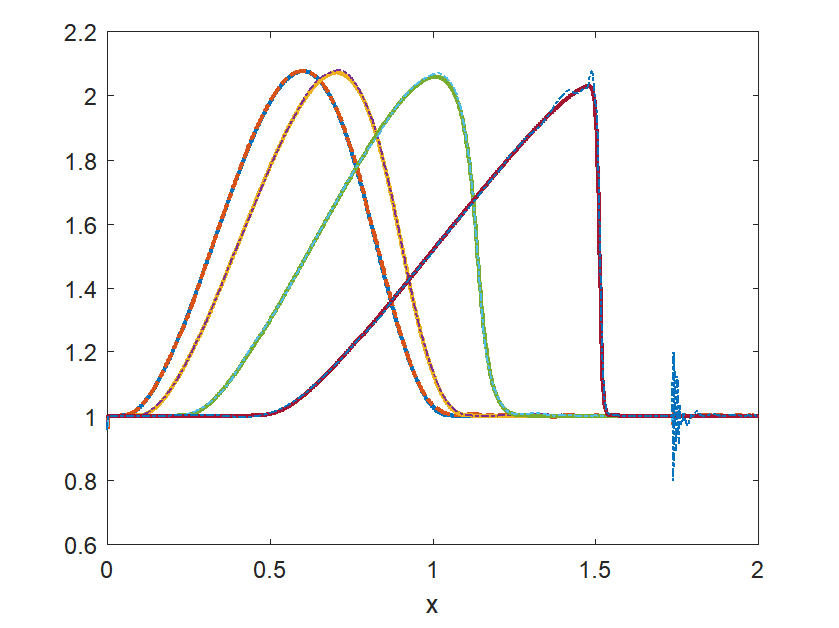}
\caption{Left: the FOM solutions (solid lines), the hybrid ROM solutions (dash lines), and the POD-FFNN-e-decoder predictions (dotted lines) at time instances 0.2s, 0.4s, 1s and 2s and at at $\mu^*=1.08$. Right: the FOM solutions (solid lines) and the CAE-FFNN predictions (dash-dot lines) at time instances 0.2s, 0.4s, 1s and 2s and at at $\mu^*=1.08$}
\label{fig:solution_comp2}
\end{figure}

The trained G-ROM, hybrid ROM, POD-FFNN-e-decoder, and CAE-FFNN are used to predict the FOM solution at four testing parameter samples: $[0.92, 0.98, 1.02, 1.08]$ that are different from the training ones. The FOM solutions produced by the hybrid ROM as well as the POD-FFNN-e-decoder, at $\mu=1.08$, are plotted in Fig.~\ref{fig:solution_comp1}, where the difference of the predicted solution from the reference FOM solution is invisible. The FOM solutions at different time instances, those predicted by the hybrid ROM and the POD-FFNN-e-decoder are presented in the left plot of Fig.~\ref{fig:solution_comp2}. We can see the the POD-FFNN-e-decoder is less accurate at the last time step. The right plot of Fig.~\ref{fig:solution_comp2} compares the FOM solutions with those predicted by the CAE-FFNN, where the big error of CAE-FFNN at the last time step is obvious. 

The errors of the approximate solutions computed from the RBM ROM ($\tilde r=44$), the errors of the solutions predicted by the CAE-FFNN, as well as the errors of the solutions computed from the hybrid ROM, and those from its non-intrusive variant, POD-FFNN-e-decoder, are listed in Table~\ref{tab:errors}, respectively. Although with 1000 epochs, 1dCNNResi of the hybrid ROM reaches a loss value of O($10^{-4}$). Whereas, the CAE-FFNN only attains at O($10^{-3}$) after 1000 epochs.  After 1000 epochs, 1dCNNResi, acting as the e-decoder of the POD-FFNN-e-decoder also reaches a loss value of O($10^{-3}$). 

It can be seen that the proposed hybrid ROM performs well for different testing parameters, and is consistently one-order of magnitude more accurate than the CAE-FFNN surrogate, and is also more accurate than its non-intrusive variant.This can be attributed to the physics-informed nature of the hybrid-ROM, where the input to the NN component is the solution to the model component, which reflects the physics of the FOM. In other words, the latent space of the NN component is obtained from a structure-preserving physical model, whereas the latent space of either the CAE or the POD-FFNN-e-decoder is generated only based on data. It can be further observed that 
POD-FFNN-e-decoder is also more accurate than the CAE-FFNN, though both are data-driven. This is mainly because that POD-FFNN-e-decoder learns the error and uses the learned error to further calibrate approximate solution, while CAE-FFNN directly learns the solution. The step of calibration plays a key role for improved the accuracy. The RBM ROM also gives accurate results with small errors comparable to those of the G-ROM+1dCNNResi, but with much larger reduced dimension $\tilde r=44$. 

\begin{table}[!h]
\centering
\caption{Errors: the RBM ROM, the CAE-FFNN, the POD-FFNN-e-decoder, and the hybrid ROM (G-ROM+1dCNNResi).}
\label{tab:errors}
\begin{tabular}{|l|l|l|l|l|}
\hline
Method & $\mu=0.92$ & $\mu=0.98$ & $\mu=1.02$ & $\mu=1.08$   \\ \hline
RBM ROM& $1.2\times 10^{-3}$ & $9\times 10^{-4}$ & $6\times 10^{-4}$ & $2\times 10^{-4}$   \\ \hline
CAE-FFNN &$7.1\times 10^{-3}$ & $8.1\times 10^{-3}$ & $8.2 \times 10^{-3}$ & $8.2 \times 10^{-3}$ \\ \hline
POD-FFNN-e-decoder & $1.6\times 10^{-3}$ & $1.7\times 10^{-3}$ & $1.6\times 10^{-3}$ & $2.1\times 10^{-3}$ \\ \hline 
G-ROM+1dCNNResi & $9.2\times 10^{-4}$ &$ 6.6\times 10^{-4}$ & $7.2\times 10^{-4}$ & $7.1\times 10^{-4}$   \\ \hline
\end{tabular}
\end{table}

In Table~\ref{tab:time}, we list the offline training time and the online prediction time of all methods. The prediction time is the wall-clock time of running the surrogate model to get approximate solutions at a given testing parameter and {\it all} $n_t=$500 time steps. The training time of the hybrid ROM: G-ROM+1dCNNResi includes the POD-greedy runtime, the EIM offline computation time, and the NN-part training time. On the other hand the training time of CAE-FFNN is only the NN training time. Training POD-FFNN-e-decoder, we need apply POD to the training data to get the reduced state. Therefore, beside the NN training time, it contains also the POD compression time. 

We observe that the offline training time of either the 1dCNNResi of the hybrid ROM or the POD-FFNN-e-decoder is less than that of CAE-FFNN, since 1dCNNResi has far less NN parameters than the CAE-FFNN, thanks to the light structure of 1dCNNResi. The POD compression time during POD-FFNN-e-decoder training is only a very small part of the total training time. All surrogate models have shorter online prediction times compared to the FOM simulation time. Note that the G-ROM from RBM is faster than all other methods. However, for the much larger model in the next section, the online performance will be different. With this example, we aim to demonstrate that the proposed methods outperform the existing CAE-FFNN method in terms of accuracy and training time, while remain comparable online prediction time.
\begin{table}[!h]
\centering
\caption{Runtime:  the RBM ROM, the CAE-FFNN, the POD-FFNN-e-decoder, and the hybrid ROM (G-ROM+1dCNNResi).}
\label{tab:time}
\begin{tabular}{|l|c|c|c|c|c|}
\hline
Method &  \multicolumn{4}{|c|}{Training time} & Prediction time   \\ \hline
           & POD-Geedy & EIM  & POD & NN  &  \\ \hline
FOM & -- & --  & --    & -- &  0.91s \\ \hline
G-ROM & 6.62s & 1.16s  & -- & --   &  0.076s \\ \hline
CAE-FFNN  & -- & --  &  -- & 39m &  0.168s   \\ \hline
POD-FFNN-e-decoder & -- & --  & 0.27s & 26m &  0.185s   \\ \hline
G-ROM+1dCNNResi & 0.82 s\footnote{This time is estimated from the POD-Greedy runtime for constructing the RBM ROM using $\tilde V \in \mathbb R^{n\times \tilde r}, \tilde r=44.$ Since we take the first $r=5$ columns from $\tilde V$ to get $V$, the runtime of POG-Greedy for getting $V$ is estimated as $1/8$ of the total 6.62 seconds in the third line}  & 1.16s  & -- & 19m   &      0.16s              \\ \hline
\end{tabular}
\end{table}

\subsection{2D inviscid Burgers' equation}
The 2D inviscid Burgers' equation is given as 
\begin{equation}
	\begin{aligned}
		\frac{\partial u_x}{\partial t} 
		+ \frac{1}{2}\frac{\partial u_x^2}{\partial x} 
		+ \frac{1}{2}\frac{\partial (u_x u_y)}{\partial y} &= 0.02 \, e^{\mu_2 x}, \\
		\frac{\partial u_y}{\partial t} 
		+ \frac{1}{2}\frac{\partial u_y^2}{\partial y} 
		+ \frac{1}{2}\frac{\partial (u_x u_y)}{\partial x} &= 0,  (x, y) \in \Omega:=[0,100]^2, t \in (0,25],
	\end{aligned}
	\end{equation}
with boundary conditions
\begin{align*}
		u_x(x, y, t; \mu) = \mu_1, && (x, y) \in \Gamma_1:=[x=0, y\in [0, 100] ], \\
                  u_y(x,y, t; \mu)=0,  && (x, y) \in \Gamma_1, \\
		u_x(x,y, t; \mu)=u_y(x,y, t; \mu)=0,  && (x, y) \in \Gamma_2:=[y=0, x\in [0, 100] ], \\
                   \frac{\partial u_y(x,y,t; \mu)}{\partial n}=\frac{\partial u_x(x,y,t; \mu)}{\partial n}=0, &&  (x, y) \in \Gamma \setminus \Gamma_1\cup  \Gamma_2,
 \end{align*}
where $\Gamma$ is the boundary of $\Omega$. 
The initial conditions are
$$u_x(x, y, 0; \mu) = 1, u_y(x, y, 0; \mu) = 1.$$
The vector of two parameters is $\mu:=(\mu_1, \mu_2)$, where $\mu_1$ defines the parametric boundary condition, and $\mu_2$ defines the parametrized source term.  $(\mu_1, \mu_2) 
\in [4.25, 5.50] \times [0.015, 0.03]$. 

We use finite difference method to discretize the 2D Burgers' equation in space, where the backward finite difference scheme is used. The x- and y-direction of $\Omega$ is divided into 249 elements, respectively, resulting in 250 grid points in each direction. Finally, we get a FOM with $n=248 \times 248 \times 2 =123,008$, by excluding the two boundary points in each direction, respectively. However, the influence of $\mu_1$ on the  solution at other grid points remains, because of the backward finite difference scheme.

\subsubsection{Training the G-ROM, the LSPG-ROM and the model component of the hybrid ROM}
We use the following settings for the POD-Greedy algorithm to obtain the G-ROM, the LSPG-ROM and the model component of the hybrid ROM:
\begin{itemize}
  \item Tolerance for the greedy algorithm: $tol_{RB} = 2.4$  
    \item Tolerance for the SVD decomposition during POD-greedy: $tol_{svd} = 0.01$.  
  \item Tolerance for the EIM: $tol_{EIM}=0.7$.  
    \item Training parameter sample set: $\mathcal P_{train}$ is composed of tensor grid points in the parameter domain by taking 3 equidistant points from the $\mu_1$-interval and the $\mu_2$-interval, respectively.  
\end{itemize}
When we use the G-ROM to compute the error estimator $\eta(\mu)$ at each iteration of Algorithm~\ref{alg:POD-greedy}, we get a projection matrix $\tilde V \in \mathbb R^{n\times \tilde r}, \tilde r=227$, from which we can construct a RBM G-ROM with size $\tilde r=227$. 
Alternatively, we also use the LSPG-ROM to compute the error estimator and finally obtain a $\tilde V \in \mathbb R^{n\times \tilde r}, \tilde r= 287$. A RBM LSPG-ROM can then be constructed using $\tilde V$.  These two ROMs are later used for performance comparison. 

We then compute SVD of the snapshot matrix $X_s$ assembled by all the snapshots corresponding to the POD-greedy-selected parameter samples for getting the above G-ROM. The first $r=10$ dominant left singular vectors are used to form the matrix $V$ in~(\ref{eq:rom_EI}) and get the model component of the hybrid model. The model component as a $10$-dimensional G-ROM then can be simulated in the time domain using the backward Euler method combined with the Newton's method.

In order to have a direct comparison with the methods in~\cite{morBarFM2023,morParTHetal26}, we further construct two LSPG-ROMs with $r=10$ and $r=150$, respectively. To this end, we implement the POD-greedy algorithm with an additional setting $r_{\mathrm{max}}=10$ or $r_{\mathrm{max}}=150$. When the number of columns in $\tilde V$ is larger than $\tilde r$, POD-greedy algorithm stops. The error estimator is computed from the LSPG-ROM built at each iteration.  The finally obtained $\tilde V$ from the POD-Greedy algorithm may not have exactly 10, or 150 columns. We further truncate the columns of $\tilde V$ to obtain the LSPG-ROMs with $r=10$ and $r=150$, respectively. We compare our proposed surrogate with these two ROMs. In~\cite{morBarFM2023,morParTHetal26}, $\tilde r=150$ is the dimension of the full linear subspace Range($\tilde V$) comprising the 10-dimensional latent space Range $(V_1)$ defining the ROM and the remaining 140-dimensional linear subspace Range ($V_2$) used for nonlinear manifold construction (see the explanations in Section~\ref{sec:NManifold}). With these two LSPG-ROMs, we know that
\begin{itemize}
 \item Since the nonlinear manifold of the reduced states is an approximation of the 140-dimensional linear subspace, the finally constructed 10-dimensional ROMs with nonlinear manifold augmentation~\cite{morBarFM2023,morParTHetal26} are no more accurate than the LSPG-ROM with $r=150$. 
 
\item Since the 10-dimensional ROM (r=10) with nonlinear manifold augmentation~\cite{morBarFM2023,morParTHetal26} includes not only the nonlinear operators present in LSPG-ROMs but also additional nonlinear functions of the reduced states arising from the nonlinear manifold construction, it is significantly more complex than a standard LSPG-ROM with the same dimension. Consequently, the online runtime of the 10-dimensional LSPG-ROM is expected to be substantially shorter than that of the 10-dimensional ROMs with nonlinear manifold augmentation.
\end{itemize}
Finally, both ROMs can serve as reference surrogate models for comparison with our proposed surrogates and those in~\cite{morBarFM2023,morParTHetal26}. 

\subsubsection{Training the NN component of the hybrid ROM}
\label{sec:NN-2d}
For the 2D inviscid Burgers' equation, we try to explore the spatial information during deep learning. Therefore, the residual block uses Conv2d layers instead of Conv1d. It is plotted in Fig.~\ref{fig:resi_block_2d}, where the output is $\mathrm{2dResiblock}({\bf x})$. The upsampling block shown in Fig.~\ref{fig:up_block_2d} contains a ConvTranspose2d layer and two residual blocks 2dResiblock. The final 2dCNNResi in Fig.~\ref{fig:Unet} is composed of 2 2dResiblocks and 3 upsampling blocks. For the 2D inviscid Burger's equation, the solution vector contains the velocity components in the x- and y-direction. If one is interested in the velocity components in both directions, the error snapshot data can be separated into 2 channels, representing the data in the two different directions, respectively. Another choice is to use two separate NNs to learn each component of the velocity, which can be implemented in parallel. One the other hand, if one is only interested in how the x-directional solution $u_x$ changes directly with respect to the two parameters, then only the snapshot data of $u_x$ are used to train the NN. This can further reduce the number of NN parameters. Here, we train the NN component using only the snapshots of $u_x$. 

The training data for the 2dCNNResi include the snapshots of the solution error contributed by the model component, and the snapshots of the reduced state vector ${\bf z}(t, \mu)$ at the training parameter samples and time instances. The training parameter samples are the same samples in $\mathcal P_{train}$, which are also the training samples used in~\cite{morParTHetal26,morBarFM2023}. On the other hand, the training data for POD-2dCNNResi are the POD coefficients $\hat {\bf e}(t, \mu)$ of the solution error, which are computed from~(\ref{eq:PODcoefficients}). Here, we take $V_0\in \mathbb R^{n\times r_e}, r_0=1024$. 

Both the FOM and the model component of the hybrid ROM are simulated using the backward Euler scheme and the Newton's method in the time domain. Then the same process for computing the error snapshots ${\bf e}(t, \mu)$ for the 1D Burgers' equation can be done for the 2D case. 

\subsubsection{Training the POD-FFNN-e-decoder}
We use the same POD-2dCNNResi in subsection~\ref{sec:NN-2d} as the e-decoder of the POD-FFNN-e-decoder. The same training parameter samples for the hybrid ROM are used to generate the training data for the POD-FFNN-e-decoder. We apply SVD to the matrix of the solution snapshots to get the snapshots of ${\bf z}(t, \mu) \in \mathbb R^r, r=10$ in Fig.~\ref{fig:non-intrusive} following~(\ref{eq:solution-PODcoefficients}). These are the training data for the FFNN. 
The training data for the e-decoder can be obtained according to~(\ref{eq:PODcoefficients}). To this end, we first compute SVD of the matrix $X_e$ in~(\ref{eq:errorsnapshots_edecoder}). Let $V_0$ be a matrix composed of the $r_0=1024$ dominant left singular vectors of $X_e$. The POD coefficient matrix ${Z}_e$ of the error snapshots $X_e$ is computed as ${Z}_e=V_0^TX_e, {Z}_e:=[{\bf z}_e(t_0, \mu_1), \ldots, {\bf z}_e(t_{n_t-1}, \mu_1) \ldots, {\bf z}_e(t_0, \mu_{n_\mu}), \ldots, {\bf z}_e(t_{n_t-1}, \mu_{n_\mu})]$, which are used as the output training data for the e-decoder. $loss_2$ is defined as the difference between ${\bf z}_e(t, \mu)$ and $\hat e(t, \mu)$, the output of the neural network. The final output of the e-decoder is computed as $\tilde {\bf e}(t, \mu)=V_0 \hat e(t, \mu)$ according to~(\ref{eq:PODcoefficients}).

\subsubsection{Prediction results}
We present the results of two hybrid-ROMs: G-ROM+2dCNNResi, G-ROM+POD-2dCNNResi and those of the POD-FFNN-e-decoder. We also compare them with the results of G-ROM and LSPG-ROM. The only difference between the G-ROM+2dCNNResi surrogate and the G-ROM+POD-2dCNNResi surrogate is the NN component, for the G-ROM+2dCNNResi surrogate, the error $\tilde {\bf e}(t, \mu)$ is learned directly by a 2dCNNResi network, whereas for the G-ROM+POD-2dCNNResi surrogate, the error $\tilde {\bf e}(t, \mu)$ is learned by POD-2dCNNResi. 

All the results are obtained at the three testing samples $${\mathcal P}_{test}=\{(5.19, 0.026), (4.75, 0.02), (4.56, 0.019)\}.$$ The same testing samples are used in~\cite{morParTHetal26,morBarFM2023}. This gives a direct comparison to the results presented in~\cite{morParTHetal26,morBarFM2023}. The 2dCNNResi of G-ROM+2dCNNResi is trained for 100 epochs, while the 2dCNNResi of  G-ROM+POD-2dCNNResi is trained for 5500 epochs. The POD-FFNN-e-decoder is trained for 10,000 epochs.

We present in Table~\ref{tab:2derrors} the errors of the following surrogate models: the G-ROM with $r=227$, a LSPG-ROM with $r=10$, a LSPG-ROM with $r=150$, a LSPG-ROM with $r=287$, the G-ROM+2dCNNResi, the G-ROM+POD-2dCNNResi and the POD-FFNN-e-decoder. 
In Table~\ref{tab:2dtime}, we list the offline training time and the online prediction time of the models in Table~\ref{tab:2derrors}. The prediction time is the wall-clock time of running the surrogate model to get approximate solutions at a given testing parameter and at {\it all} 500 time steps. The training time includes the POD-greedy runtime, the EIM offline computation time, the POD compression time if applicable, and the NN training time. 

\begin{figure}[!h]
\centering
\includegraphics[width=140mm]{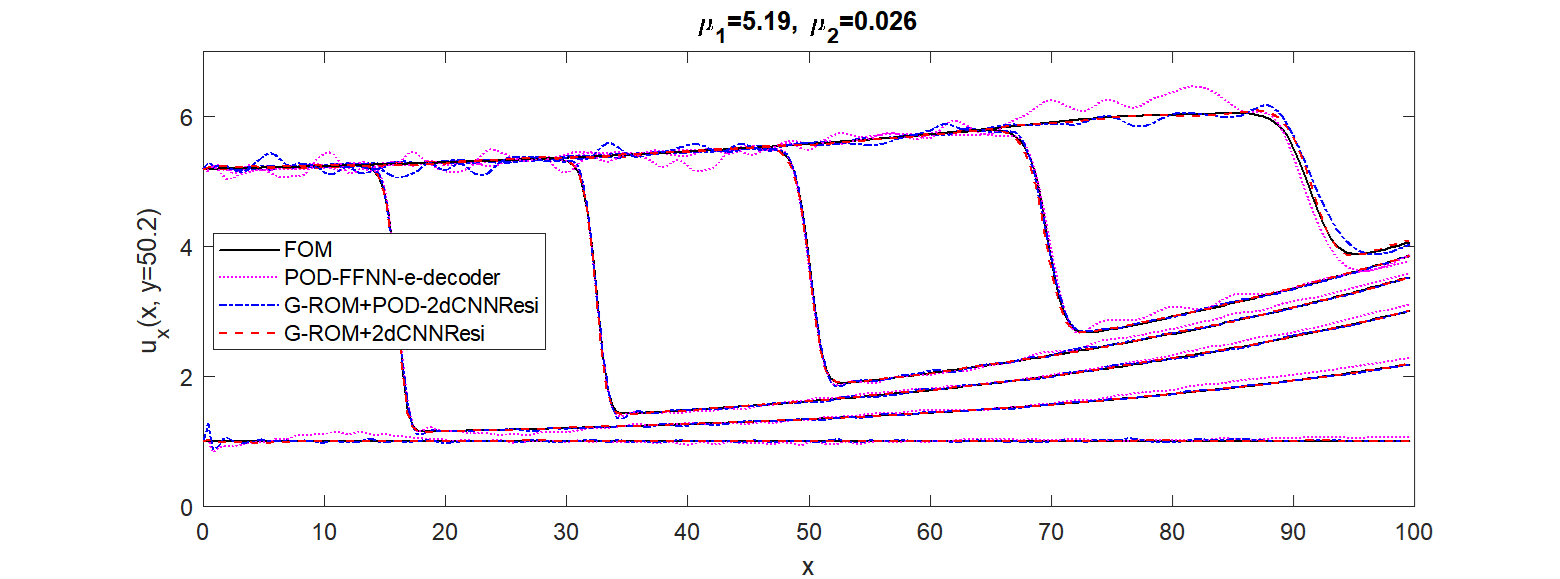} \quad
\includegraphics[width=140mm]{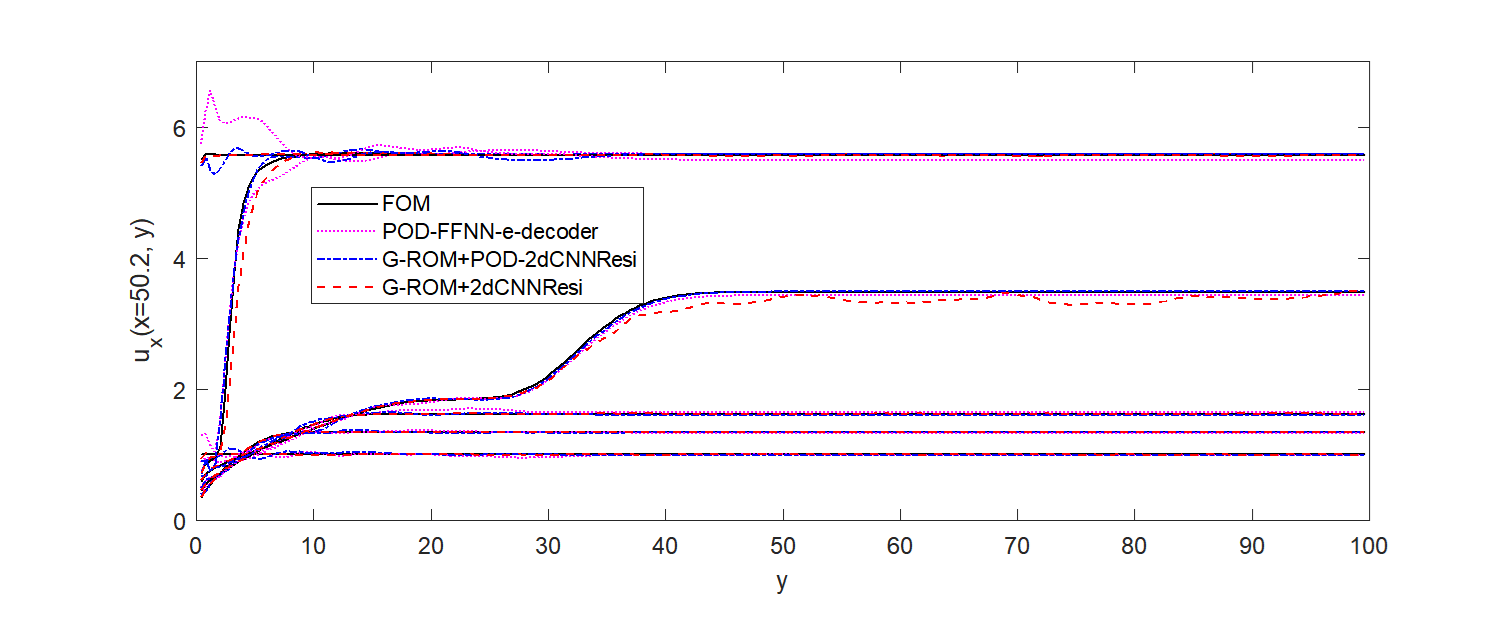}
\caption{The FOM solutions and the surrogate predictions at time instances 0, 5, 10, 15, 20, and 25.}
\label{fig:para1}
\end{figure}

\begin{figure}[!h]
\centering
\includegraphics[width=140mm]{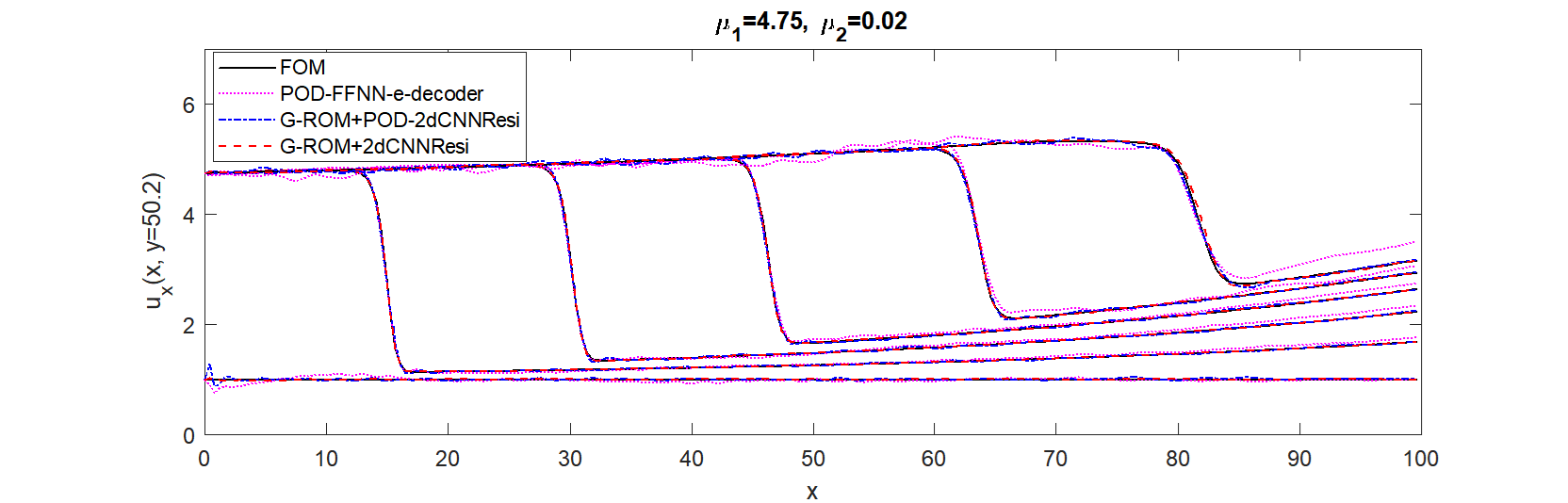} \quad
\includegraphics[width=140mm]{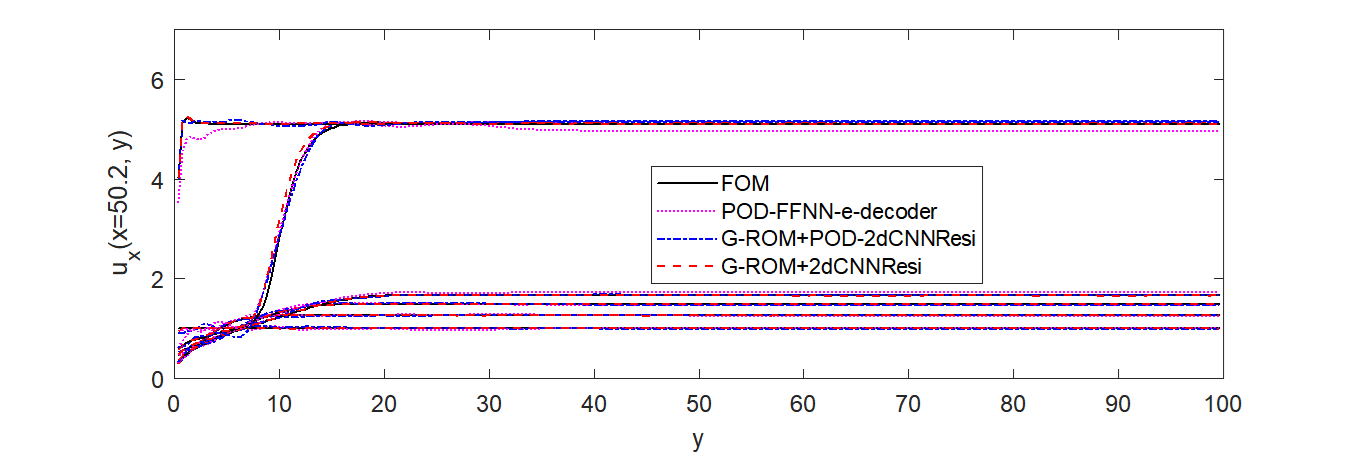}
\caption{The FOM solutions and the surrogate predictions at time instances 0, 5, 10, 15, 20, and 25.}
\label{fig:para2}
\end{figure}

\begin{figure}[!h]
\centering
\includegraphics[width=140mm]{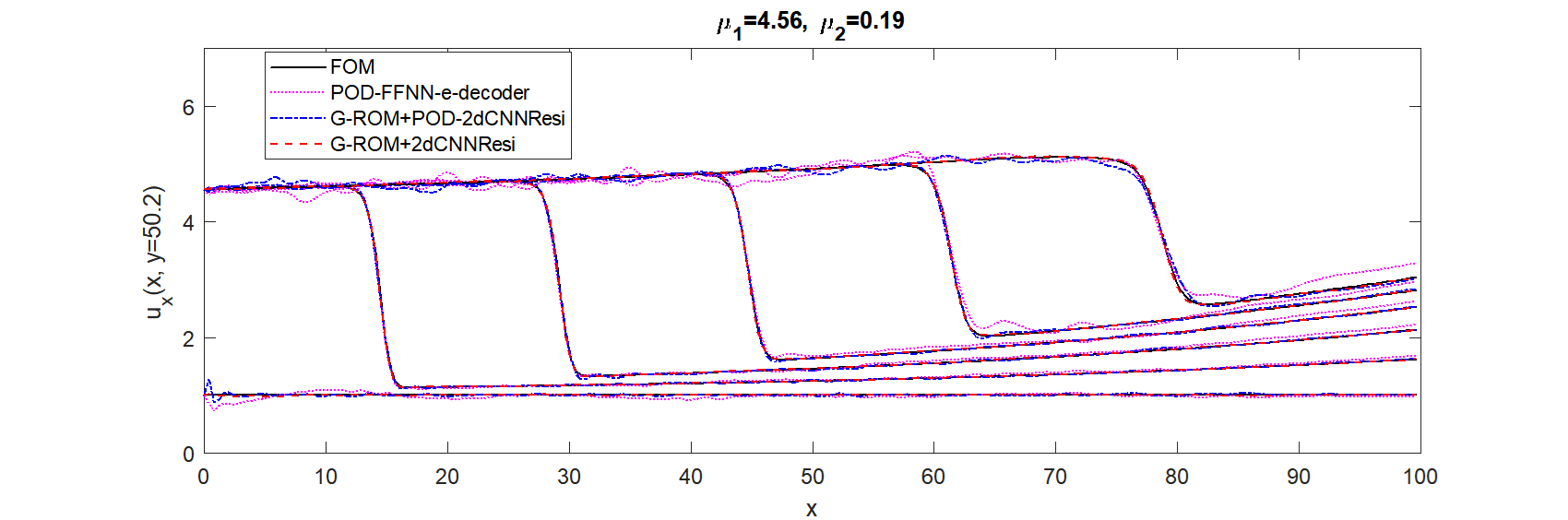} \quad
\includegraphics[width=140mm]{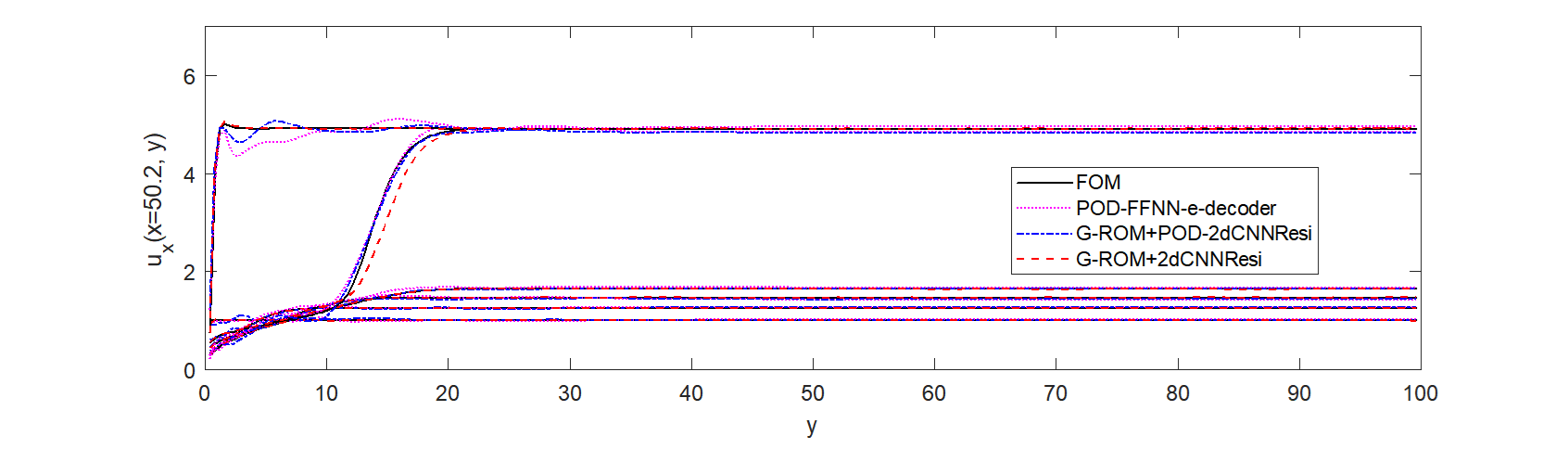}
\caption{The FOM solutions and the surrogate predictions at time instances 0, 5, 10, 15, 20, and 25.}
\label{fig:para3}
\end{figure}

Figures~\ref{fig:para1}-\ref{fig:para3} plot the predicted solutions along one space direction by fixing another space direction at the times $t=0, 5, 10, 15, 20, 25$. These again provide a direct comparison between the proposed method and the methods in~\cite{morParTHetal26,morBarFM2023}, where corresponding results of ROMs with nonlinear manifold augmentation are also shown.

By summarizing the results in both the Tables and the Figures, we obtain the following observations:
\begin{itemize}
 
\item On the one hand, the hybrid ROM, G-ROM+POD-2dCNNResi, is more accurate than the LSPG-ROM with $r=150$. This means it is definitely more accurate than the ROMs with nonlinear manifolds in~\cite{morBarFM2023,morParTHetal26}, where the total linear subspace dimension is $\tilde r=150$, the dimension of the ROM is $r=10$, and the dimension of the linear subspace Range$(V_2)$ for nonlinear manifold augmentation (see subsection~\ref{sec:NManifold}) is 140. This justifies our analysis in subsection~\ref{sec:NManifold} that the hybrid ROM breaks the barrier of linear-subspace-limited-accuracy of the ROMs with nonlinear manifold augmentation. On the other   hand, it is 89x faster than the LSPG-ROM with $r=10$ for online prediction, which implicates that it is at least 89x faster than the ROMs with nonlinear manifolds~\cite{morBarFM2023,morParTHetal26}. The online prediction time of G-ROM+2dCNNResi is longer, but still not much longer than the LSPG-ROM with $r=10$. This indicates that G-ROM+2dCNNResi might have comparable online prediction time as the ROMs in~\cite{morBarFM2023,morParTHetal26}.

\item For the G-ROM with $r=227$, it is also less accurate than the proposed surrogate models, and delivers less speed-ups compared to the G-ROM+POD-2dCNNResi and the POD-FFNN-e-decoder. 

\item The LSPG-ROM with $r=287$ is most accurate, but has no speed-up for online prediction.

\item After training with 100 epochs, the G-ROM+2dCNNResi surrogate is already accurate enough. However, it has very long training time. In contrast, the G-ROM+POD-2dCNNResi and POD-FFNN-e-decoder attain similar accuracy with much shorter training time. 

\item  From the figures, we see that G-ROM+2dCNNResi surrogate gives the most accurate and smooth prediction. POD-FFNN-e-decoder has larger oscillations at the last time instance $t=25$. Despite that, the results of G-ROM+2dCNNResi and those of G-ROM+POD-2dCNNResi show an overall improvement in accuracy as compared to those results in the Figures in~\cite{morParTHetal26,morBarFM2023}. 
\end{itemize}
In summary, the proposed two surrogates G-ROM+POD-2dCNNResi and POD-FFNN-e-decoder achieve a three-order-of-magnitude speedup in online prediction, and are up to two orders of magnitude faster than the ROMs with nonlinear manifolds~\cite{morParTHetal26,morBarFM2023}, while maintaining acceptable offline training costs and achieving improved accuracy. The non-intrusive variant POD-FFNN-e-decoder provides greater flexibility for surrogate modeling and is particularly promising for large-scale problems simulated by commercial software as black boxes. The maximal error $1.54 \times 10^{-2}$ of POD-FFNN-e-decoder is slightly smaller than the errors of the intrusive methods in~\cite{morParTHetal26,morBarFM2023}, which are $1.83 \times 10^{-2}$, $1.71 \times 10^{-2}$, and $1.66 \times 10^{-2}$, respectively. 

\begin{table}[!h]
\centering
\caption{Errors and speedups of different ROMs at testing samples of $\mu$.}
\label{tab:2derrors}
\begin{tabular}{|l|l|l|l|l|l|}
\hline
Method & $(5.19, 0.026)$ & $(4.75, 0.02)$ & $(4.56, 0.019)$ & $\epsilon_{\max}$ & Speedup \\ \hline
G-ROM (r=227) & $9.8\times 10^{-3}$ & $2.3\times 10^{-2}$ & $2.8\times 10^{-2}$ &  $2.8\times 10^{-2}$ & 16.9x \\ \hline
LSPG-ROM (r=10) & $1.4\times 10^{-1}$ & $1.5\times 10^{-1}$ & $1.4\times 10^{-1}$ &  $1.5\times 10^{-1}$  & 10.2x  \\ \hline
LSPG-ROM (r=150) & $4.0\times 10^{-3}$ & $7.1\times 10^{-3}$ & $9.4\times 10^{-3}$ &  $9.4\times 10^{-3}$ & 2x \\ \hline
LSPG-ROM (r=287) & $1.2\times 10^{-3}$ & $8.9\times 10^{-4}$ & $1.3\times 10^{-3}$ &  $1.3\times 10^{-3}$ & 1x \\ \hline
G-ROM+2dCNNResi &$1 \times 10^{-2}$  & $5\times 10^{-3}$  & $8.2\times 10^{-3}$ & $1\times 10^{-2}$  & 7.2x \\ \hline
G-ROM+POD-2dCNNResi & $ 7.2\times 10^{-3}$ & $ 6.8\times 10^{-3}$ & $ 7.8\times 10^{-3}$ & ${\bf 7.8\times 10^{-3}}$  & {\bf 248x} \\ \hline
POD-FFNN-e-decoder& $1.48 \times 10^{-2}$ & $1.51 \times 10^{-2}$ & $1.54 \times 10^{-2}$ &$1.54 \times 10^{-2}$ & {\bf 1945x} \\ \hline
\end{tabular}
\end{table}
\begin{table}[!h]
\centering
\caption{Runtime of the proposed methods}
\label{tab:2dtime}
\begin{tabular}{|l|c|c|c|c|}
\hline
Method &  \multicolumn{3}{|c|}{Training time} & Prediction time  \\ \hline
           & POD-Geedy & EIM  & (POD+) NN    &   \\ \hline
FOM & --  & --  &   --  & 233.5 s               \\ \hline
G-ROM (r=227) & 244s  & 3.6h  & --    &    13.7s          \\ \hline
LSPG-ROM (r=10) & 223s  & 3.6h  & --    &    23s         \\ \hline
LSPG-ROM (r=150) & 10.9m  & 3.6h  & --    &    112s       \\ \hline
LSPG-ROM (r=287) & 22.6m  & 3.6h  & --    &    234.8s        \\ \hline
G-ROM+2dCNNResi & 244s  & 3.6h  & 68h    &    32.4s               \\ \hline
G-ROM+POD-2dCNNResi  & 244s  & 3.6h    &  8h&  {\bf 0.96s}  \\ \hline
POD-FFNN-e-decoder  & --  & --    &  12.6h&  {\bf 0.12s}  \\ \hline
\end{tabular}
\end{table}

\section{Conclusions}
In this work, we propose two surrogate modeling methods that are efficient for convection-dominated parametric time-dependent problems. Comparison with the state-of-the-art methods shows that the proposed surrogates produce more accurate predictions and enable much faster online prediction for large-scale problems. Thanks to the physics-informed latent space, the hybrid ROM is more accurate than its non-intrusive variant. Nevertheless, the non-intrusive variant provides greater flexibility and is particularly well suited for complex problems where discretized operators are not easily accessible. Future work may include reducing the offline training costs by developing more efficient neural networks structures, further improving the accuracy of the non-intrusive surrogate, and investigating techniques for optimally determining the latent space dimension.

\newpage
\bibliographystyle{abbrv}
\bibliography{refs_updated2}

\newpage 
\appendix
\section{Details of the NN structure for the examples}
\label{sec:appendix}

\begin{figure}[!h]
\centering
\includegraphics[width=100mm]{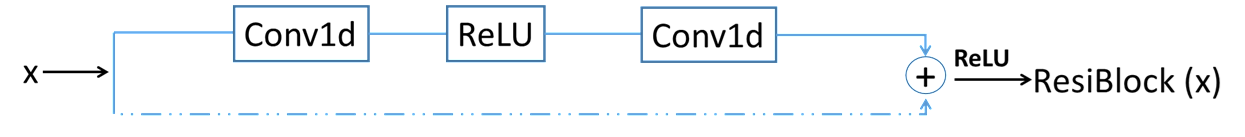}
\caption{Residual block for the 1D inviscid Burgers' equation}
\label{fig:resi_block}
\end{figure}
\ \\

\begin{figure}[!h]
\centering
\includegraphics[width=100mm]{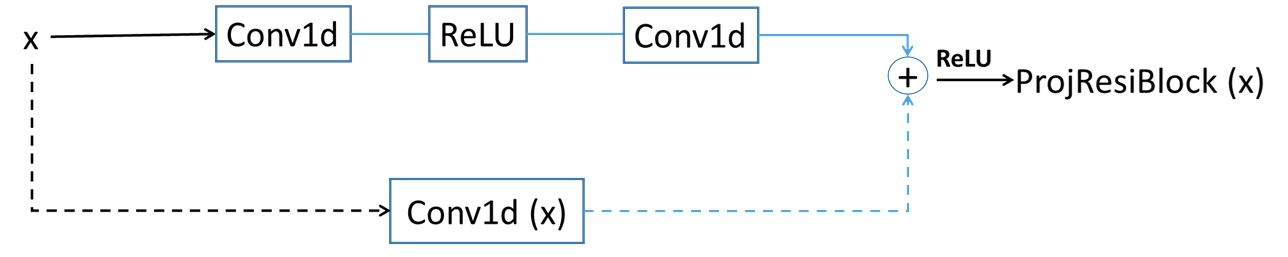}
\caption{Residual block used in the upsampling block for the 1D inviscid Burgers' equation}
\label{fig:Projresi_block}
\end{figure}
\ \\

\begin{figure}[!h]
\centering
\includegraphics[width=80mm]{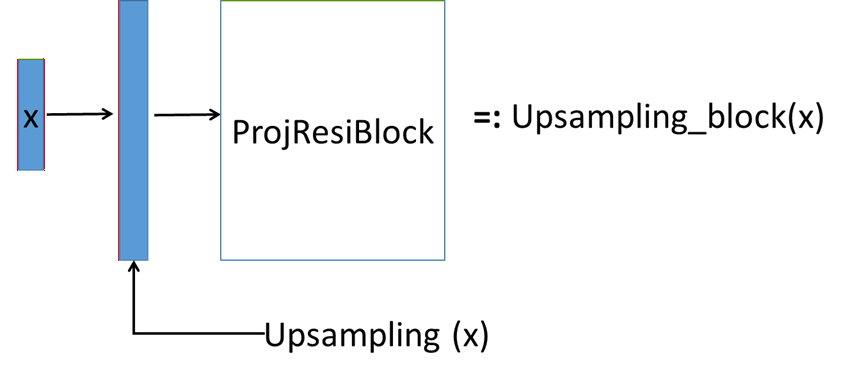}
\caption{Upsampling block for the 1D inviscid Burgers' equation}
\label{fig:up_block}
\end{figure}
\ \\

\begin{figure}[!h]
\centering
\includegraphics[width=100mm]{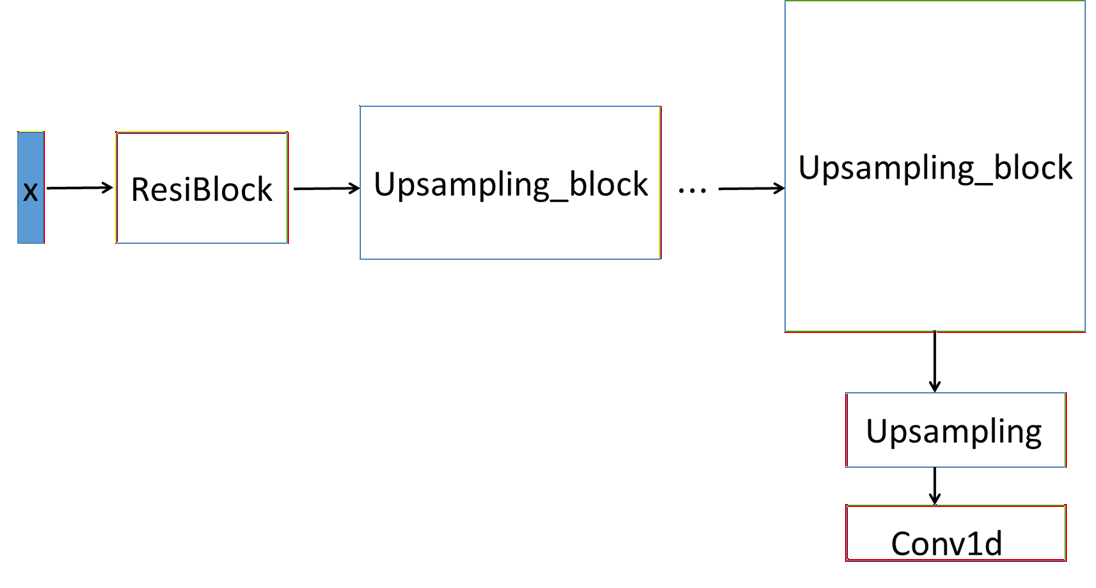}
\caption{The 1dCNNResi for the 1D inviscid Burgers' equation}
\label{fig:CNNresi}
\end{figure}
\ \\

\begin{figure}[!h]
\centering
\includegraphics[width=100mm]{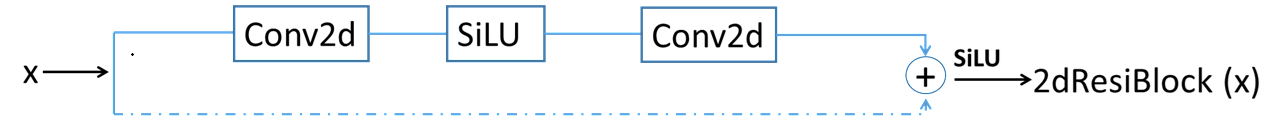}
\caption{Residual block for the 2D inviscid Burgers' equation}
\label{fig:resi_block_2d}
\end{figure}
\ \\

\begin{figure}[!h]
\centering
\includegraphics[width=100mm]{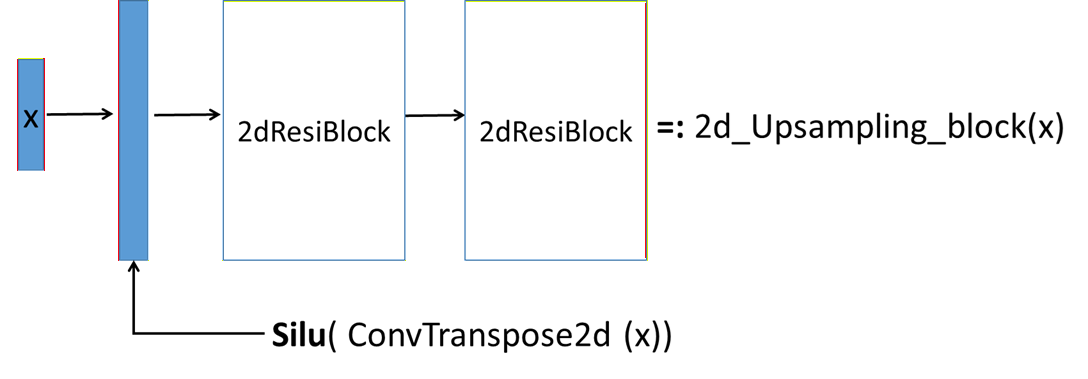}
\caption{Upsampling block used in the upsampling block for the 2D inviscid Burgers' equation}
\label{fig:up_block_2d}
\end{figure}
\ \\
\ \\

\begin{figure}[!h]
\centering
\includegraphics[width=100mm]{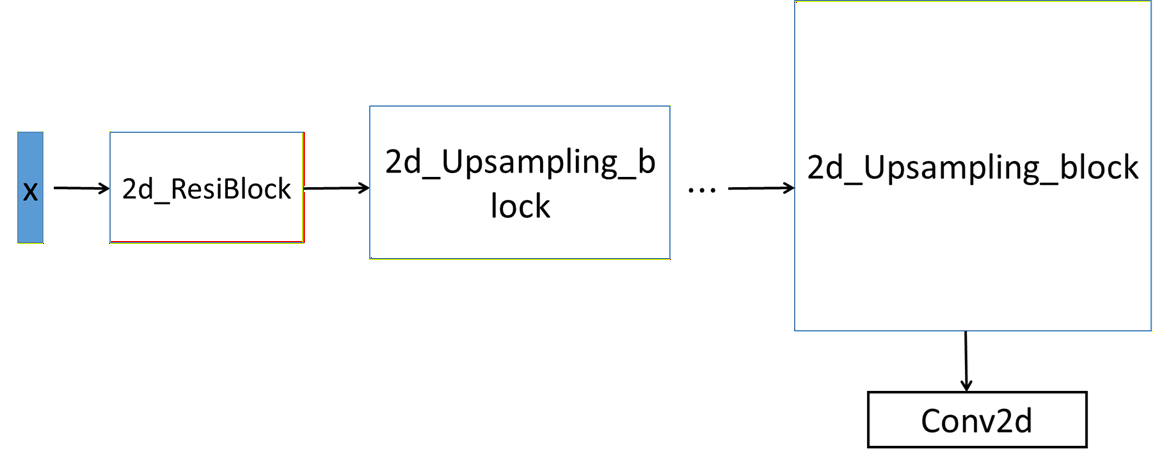}
\caption{The 2dCNNResi for the 2D inviscid Burgers' equation}
\label{fig:Unet}
\end{figure}

\end{document}